\documentclass[12pt]{amsart}

\usepackage{amsfonts, amssymb, curves}
\newcommand{\R}{\mathbb{R}}
\newcommand{\N}{\mathbb{N}}
\newcommand{\Q}{\mathbb{Q}}
\newcommand{\Z}{\mathbb{Z}}
\newcommand{\cT}{\mathbb{T}}
\newcommand{\T}{\mathcal{T}}

\newcommand{\A}{\mathcal{A}}
\newcommand{\B}{\mathcal{B}}
\newcommand{\C}{\mathcal{C}}
\newcommand{\cL}{\mathcal{L}}

\newcommand{\cP}{\mathcal{P}}

\newcommand{\V}{\mathcal{V}}
\newcommand{\F}{\mathcal{F}}
\newcommand{\G}{\mathcal{G}}

\newcommand{\vp}{\varphi}
\newcommand{\tv}{\tilde{\varphi}}
\newcommand{\tA}{\tilde{\A}}
\newcommand{\tT}{\tilde{T}}

\newcommand{\hf}{\hat{f}}
\newcommand{\td}{\tilde}

\newcommand{\larr}{\left( \begin{array}{c}}
\newcommand{\rarr}{\end{array} \right) }

\newcommand{\lsqarr}{\left[ \begin{array}{c}} 
\newcommand{\rsqarr}{\end{array} \right]}

\newcommand{\arrow}{\rightarrow}
\newcommand{\inv}{\varprojlim}
\newcommand{\dir}{\varinjlim}

\begin{document}

\newtheorem{theorem}{Theorem}[section]

\newtheorem{corollary}[theorem]{Corollary}
\newtheorem{lemma}[theorem]{Lemma} \newtheorem{prop}[theorem]{Proposition}
\newtheorem{example}[theorem]{Example}

\title{Proximality In Pisot tiling
spaces}

\author{Marcy Barge and Beverly Diamond}

 \begin{abstract} A substitution $\vp$ is strong Pisot if its abelianization
matrix is non-singular and all eigenvalues except the Perron-Frobenius eigenvalue
have modulus less than one.   For strong Pisot $\vp$ that satisfies a no cycle
condition and for which the translation flow on the tiling space $\T_\vp$ has
pure discrete spectrum, we describe the collection
$\T^P_\vp$ of pairs of proximal tilings in
$\T_\vp$ in a natural way as a substitution tiling space.  We show
that if $\psi$ is another such substitution, then $\T_\vp $ and $\T_\psi$ are
homeomorphic if and only if $\T^P_\vp$ and $\T^P_\psi$ are homeomorphic.  We make
use of this invariant to distinguish tiling spaces for which other known
invariants are ineffective. In addition, we show that for strong Pisot
substitutions, pure discrete spectrum of the flow on the associated tiling space
is equivalent to proximality being a closed relation on the tiling space.  
\end{abstract} 

 \maketitle

\begin{center}September 1, 2005

\end{center}

\section{Introduction}

Substitutions, and the tiling spaces associated with them, are of fundamental
interest in recent investigations in number theory (numeration systems in various
bases, diophantine approximation), physics (modelling quasicrystals), and
dynamical systems (coding hyperbolic attractors, generating Markov partitions)
(see
\cite{book} for an excellent introduction to these topics).

The subtle recurrence properties of words generated by a substitution are
expressed in the intricate details of the topology of the associated tiling
space.  Of particular interest are those substitutions whose abelianizations have
a dominant eigenvalue that is a Pisot number, largely because of the connection
with pure discrete spectrum of the tiling flow and its consequences.  There are
several algorithms that verify pure discrete spectrum, some of which apply
generally (\cite{anne},
\cite{mjarek}, \cite{sirsol}), and others of which apply  in special cases
(\cite{holsol}, \cite{ars}, \cite{sid}, \cite{bd2}).  For strong Pisot
substitutions (see \S \ref{sec-notation}), in all known examples, the translation
flow on the associated tiling space has pure discrete spectrum.  
The
conjecture that strong Pisot substitutions always produce tiling flows with pure
discrete spectrum is variously known as:  the Pure Discrete Spectrum Conjecture,
the Geometric Coincidence Conjecture, the Super Coincidence Conjecture, and,
simply, the Pisot Conjecture.  (See \cite{bertheanne} for a survey of the progress
on this conjecture, and \S \ref{sec-prox} for a statement of the Geometric
Coincidence Condition.)   

Given a strong Pisot substitution $\vp$ whose tiling flow has
pure discrete spectrum,   the tiling dynamics on $\T_\vp$ are measure
theoretically and almost topologically conjugate with Anosov dynamics on a torus
or solenoid of the appropriate dimension.  This almost conjugacy is simply the
quotient map that glues together proximal points in the tiling space $\T_\vp$.
That is, for such substitutions, the
gross topology of the tiling space is that of a torus or solenoid; the intricate
details lie in the structure of the collection of proximal points.  

In this paper we isolate this proximal structure, symbolically encode it, and
demonstrate its use as a distinguishing invariant for the topological type of
Pisot tiling spaces.  We also prove that, for strong Pisot substitutions, pure
discrete spectrum for the tiling flow is equivalent to proximality being a closed
relation.

Given a primitive and aperiodic substitution $\vp$ with associated tiling space
$\T_\vp$, tilings $T, T' \in
\T_\vp$ are asymptotic provided $dist(T - t, T'-t) \arrow 0$ as $t \arrow \infty$
or $t \arrow -\infty$.  In any such tiling space, there is a finite positive
number of arc components (i.e., composants), all of whose tilings are asymptotic
to the tilings in some other arc component.  Such arc components are called the
asymptotic composants of the tiling space, and any homeomorphism from one tiling
space to another must take asymptotic composants to asymptotic composants.  In
\cite{bd1}, we exploited this fact to develop a complete topological invariant
for 1-dimensional substitution  tiling spaces. Unfortunately, this invariant is
difficult to use.

More computable (but far from complete) invariants have recently  emerged
(\cite{toke},
\cite{mmegan}, \cite{mrichard}).  These are all cohomological in nature and
depend, in one way or another, on the interaction between some relative cocycles
associated with asymptotic composants and the cocycles of the space itself.  The
approach of this paper is to consider the less restrictive notion of proximal
composants (tilings $T, T' \in \T_\vp$ are proximal if $\inf_t dist (T - t, T' -
t) = 0$).  This typically provides a much richer collection of composants to
consider. 

For general tiling spaces, it is not clear how to formulate proximality as
a purely topological concept (that is, not tied to a particular flow). 
For example, in the case that $\vp$ is primitive, for any two tilings $T, T' \in
\T_\vp$, one can always find reparameterizations $\alpha, \alpha'$ of $\R$ so
that $\inf_{t} dist (T-\alpha(t),T'-\alpha'(t) ) = 0$.  
 However, in the case that tiling spaces $\T_\vp$
and $\T_\psi$ are strong Pisot (with pure discrete spectrum tiling flows and a no
cycle condition on periodic words), we prove that proximality is topological in
the following sense: If $h$ is a homeomorphism of $\T_\vp$ with $\T_\psi$, and
$\C,
\C'$ are composants of $\T_\vp$ containing proximal tilings, then $h(\C)$ and
$h(\C')$ are composants of $\T_\psi$ that contain proximal tilings (Theorem
\ref{homprox}).  

Our approach is to link
proximality under the tiling flow to a symbolic analog that holds not only for
tilings expressed in the language of the given substitution $\vp$ but also in the
language associated with any substitution obtained from $\vp$ by certain cuttings
and rewritings.  Once this is accomplished, proximality becomes topological,
since any homeomorphism between $\T_\vp$ and $\T_\psi$ can be isotoped to a
homeomorphism that has a symbolic interpretation in terms of substitutions $\tv$
and $\tilde{\psi}$ obtained from $\vp$ and $\psi$ by suitable cutting and
rewriting.

To achieve this symbolic linking, it is necessary to begin with a
substitution $\vp$ that is strong Pisot and has a tiling flow with pure discrete
spectrum.  For this linkage to persist through the cutting and rewriting
that yields $\tv$, we will need to know that $\tv$ is still (weakly) Pisot--the no
cycle condition on $\vp$, the subject of \S \ref{sec-cohomol}, guarantees this.

Moreover, we  are
able to precisely identify the proximal tilings in $\T_\vp$ by showing that
$$\T^P_\vp =  \{ (T, T') : T, T' \textrm{ are proximal in } \T_\vp\}$$ is itself a
substitution tiling space, with an algorithmically identifiable underlying
substitution $\vp_{EBP}$.  We obtain: 

\

{\parindent=0pt
{\bf Theorem \ref{homess}} {\em Suppose that $\vp$ and $\psi$ are strong
Pisot,  satisfy  the no  cycle condition on periodic words, and have tiling flows
with pure discrete spectrum.  Then
$\T_{\vp}$ and $
\T_{\psi}
 $ are homeomorphic if and only if $ \T_{\vp_{EBP}} $ and $ \T_{\psi_{EBP}}  $
are homeomorphic.}}

\

In \S \ref{sec-prox}, we give examples of substitutions $\vp$ and $\psi$ for
which the additional structure in the proximal tiling spaces allows us to deduce
that $\T_\vp$ and $\T_\psi$ are not homeomorphic.  (In these examples, the known
cohomological invariants do not allow one to distinguish the two tiling spaces). 
Section \ref{sec-notation} contains definitions and background material; \S
\ref{sec-cohomol} is devoted to the no cycle condition and is technical in
nature.  The main results appear in \S\S \ref{sec-prox} and \ref{sec-closed}.

\section{Notation and terminology} \label{sec-notation}

We introduce some of the notation and terminology necessary for the paper.  

Let $\A = \{1, 2, \ldots , card(\A)\}$ and $\B = \{1, 2, \ldots ,
card(\B)\}$  be  finite alphabets;   $\A^*$ will denote the collection of
finite nonempty words with letters in $\A$.   Given a map $\tau : \A \arrow
\B^*$, there is an associated transition matrix $A_\tau = (a_{ij})_{i\in
\B,  j \in \A}$ in which
$a_{ij}
$ is the number of occurrences of $i$ in the word $\tau(j)$; $A_\tau$ is called
the {\em abelianization} of the map $\tau$.  A map
$\tau :
\A
\arrow \B^*$ extends naturally to $\tau: \A^* \arrow \B^*$.

 A {\em substitution} is a map $\varphi : \A \arrow \A^*$; $\varphi$ is {\em
primitive} if
$\varphi^n(i)$ contains
$j$ for all
$i, j \in \A$ and sufficiently large $n$. Equivalently,  $\varphi$ is
primitive if and only if the matrix $A_\varphi$ is aperiodic, in which case
$A_\varphi$ has a simple  eigenvalue $\lambda_\varphi$ larger in modulus than its
remaining eigenvalues called the Perron-Frobenius  eigenvalue of $A_\varphi$ (and
$\varphi$).

 A word $w$ is {\em allowed for $\varphi$} if
and only if for each finite subword (i.e., factor)
$w'$ of
$w$, there are $i \in \A$ and $n \in \N$ such that $w'$ is a subword
 of
$\varphi^n(i)$; the {\em language of} $\vp$,  $ \cL_\vp $, is the set of 
finite allowed words for $ \vp$. Let
$W_\varphi$ denote the set of  allowed bi-infinite words for
$\varphi$. We identify the $0^{th}$
coordinate in a bi-infinite word
$w$ by either an indexing, as in $w = \ldots w_{-1} w_0w_1 \ldots$, or by use of a
decimal point (or both). Let $\sigma: W_\vp \arrow W_\vp$ denote the {\em shift
map}:  $$ \sigma( \ldots w_{-1} . w_0 w_1 \ldots) :=  \sigma( \ldots w_{-1}  w_0
.  w_1 \ldots) .$$ For
$w
\in W_\varphi$, the {\em shift class of
$w$} is the equivalence class of bi-infinite words:
$$[w] := 
\{w' \in W_\varphi : w' \textrm{ is in the shift orbit of  } w\} 
.$$ The
substitution 
$\varphi : \A \arrow \A^*$ extends to $\varphi: W_\varphi \arrow W_\varphi$ 
 where
$$\varphi(\ldots w_{-1}w_0w_1 \ldots) := \ldots
\varphi(w_{-1})\ {\bf .} \ \varphi(w_0)\varphi(w_1)
\ldots$$ as well as to a map on equivalence classes:  $$\varphi([w]) :=
[\varphi(w)].$$ 
The word 
$w$ is {\em  periodic} for $\varphi$, or $\varphi$-periodic,  if for some $m
\in
\N$, 
$$\varphi^m(w) =
\ldots
\varphi^m(w_{-1}) {\bf .} \varphi^m(w_0)\varphi^m(w_1)
\ldots = \ldots w_{-1} {\bf .} w_0w_1 \ldots.$$ 

 Each primitive
substitution $\varphi$ has at least one allowed $\varphi$-periodic bi-infinite word
which is necessarily uniformly recurrent under the shift.   A substitution $\varphi$ with precisely one periodic, hence
fixed, bi-infinite word  is called {\em proper}; $\vp$ is proper if and
only if there are
$b, e
\in
\A$ such that for all sufficiently large  $k$ 
 and all
$i \in \A$,
$\varphi^k(i) = b \ldots e$.

A primitive
substitution
$\varphi$ is {\em aperiodic} if at least one (equivalently, each)
$\varphi$-periodic  bi-infinite word is not periodic under the natural shift map,
in which case $(W_\varphi, \sigma)$ is an infinite minimal
dynamical system.  If $\vp$ is aperiodic, then the map $\varphi: W_\varphi
\arrow W_\varphi$ is one-to-one (\cite{mosse}). If
$\varphi$ is periodic (that is, primitive and not aperiodic), then
$W_\varphi$ is finite.

The  substitution $\varphi$ is {\em weak Pisot} if $\vp$ is primitive, aperiodic,
and all eigenvalues  of $A = A_\varphi$ other than the Perron-Frobenius
eigenvalue have modulus strictly less than 1; $\vp $ is {\em strong Pisot} if
 any non-dominant eigenvalue for $\vp$
has modulus strictly between 0 and 1.   If
$\varphi$ is strong Pisot, then
$\varphi$ is necessarily primitive and aperiodic (\cite{book} and \cite{hz}).  If
$\vp$ is weak Pisot,   then the (hyperbolic) linear map on
$\R^d$ defined by the matrix
$A$ has stable space $E^s$ of dimension $d-1$ and unstable space $E^u$ of dimension
1 spanned by the  positive right  Perron-Frobenius eigenvector $\omega_R$
for $A$. Also, if $\vp$ is strong Pisot, neither $E^s$ nor
$E^u$ contain elements of the integer lattice  other than the origin.  (See Chapter
1  of \cite{book}, for instance.)


 Given a primitive substitution  $\varphi : \A \arrow \A^*$ with $card(\A) =
d \geq 2$, let
$\omega_{L, \vp} := \omega_L = (\omega_1, \ldots, \omega_d)$ and
$\omega_{R, \vp} := \omega_R
$ be positive left and right eigenvectors (respectively) for the
Perron-Frobenius eigenvalue, 
$ \lambda = \lambda_\vp$, of $A$.
 The intervals $P_i=[0,\omega_i]$,
$i=1,\ldots,d$, are called {\em prototiles for $\vp$} (consider
$P_i$ to be distinct from $P_j$ for $i\ne j$ even if $\omega_i=\omega_j$).   A
{\em tiling} $T$ of $\mathbb{R}$ by the prototiles for
$\vp$ is a collection $T=
\{T_i\}_{i=-\infty}^{\infty}$ of tiles $T_i$ for which 
$\bigcup_{i=-\infty}^{\infty}T_i=\mathbb{R}$,  each $T_i$ is a translate of
some $P_j$ (in which case we say $T_i$ is {\em of type $j$}), and $T_i\cap
T_{i+1}$ is a singleton for each
$i$.  Generally we assume that the indexing is such that $0 \in T_0 \setminus
T_{1}$.  

There are occasions in this paper 
when we wish to define tilings by prototiles for $\vp$ but $\vp$ is not
primitive.  In each such case, the matrix for $\vp$ will have a unique
Perron-Frobenius eigenvector, so  prototiles and tilings will be well-defined.

If $\varphi(i) = i_1 i_2 \ldots i_{k(i)}$, then
$\lambda \omega_i = \sum_{j=1}^{k(i)}\omega_{i_j}$. Thus
$|\lambda P_i| = \sum_{j=1}^{k(i)}|P_{i_j}|$, and  
$\lambda P_i$ is tiled by $\{T_j\}_{j=1}^{{k(i)}}$,
where $T_j= P_{i_j} + \sum_{k=1}^{j-1}\omega_{i_k}$. This process is called
inflation and substitution and extends to a map $\Phi$ taking a tiling
$T=\{T_i\}_{i=-\infty}^{\infty}$ of $\mathbb{R}$ by prototiles to a new tiling,
$\Phi(T)$, of $\mathbb{R}$ by  prototiles defined by inflating,
substituting,  and suitably translating each $T_i$.  More precisely, for $w
= w_1 \ldots w_n \in
\A^*$, define
$$\cP_w + t = \{P_{w_1} + t, P_{w_2} + t + |P_{w_1}|, \ldots, P_{w_n} + t +
\Sigma_{i < n}|P_{w_i}|\}.$$  Then $\Phi(P_i + t) = \cP_{\varphi(i)} +
\lambda t$ and $\Phi(\{P_{k_i} + t_i\}_{i \in \Z}) = \cup_{i
\in \Z}(\cP_{\varphi(k_i)} +
\lambda t_i)$.  

There is a natural topology on the collection $\Sigma_{\varphi}$ of all tilings
of $\mathbb{R}$ by prototiles ($\{T_i\}_{i=-\infty}^{\infty}$ and
$\{T_i^{\prime}\}_{i=-\infty}^{\infty}$ are ``close" if there is an $\epsilon$ near
$0$ so that 
$\{T_i\}_{i=-\infty}^{\infty}$ and $\{T_i^{\prime} + \epsilon\}_{i=-\infty}^{\infty}$
are identical in a large neighborhood of $0$ (see
\cite{anderson} for details)).  The space
$\Sigma_{\varphi}$ is compact and  metrizable  with this topology and
$\Phi:\Sigma_{\varphi}\to\Sigma_{\varphi}$ is continuous.  
Given $T = \{T_i\}_{i = - \infty}^\infty \in \Sigma_\vp$, let
$\underline{w}(T) = \ldots w_{-1}w_0w_1 \ldots$ denote the bi-infinite word
with $w_i = j$ if and only if $T_i$ is of type $j$.   
 The  {\em
tiling space associated with}
$\varphi$, $\T_{\varphi}$, is defined as $$\T_\vp = \{ T :
\underline{w}(T) \textrm{ is allowed for } \vp\}.$$ 

 There is a natural
flow (translation) on 
$\Sigma_{\varphi}$ defined by
$(\{T_i\}_{i=-\infty}^{\infty}, t)\mapsto\{T_i - t\}_{i=-\infty}^{\infty}$. If
$\vp$ is primitive and aperiodic,
$\Phi:\T_{\varphi}\to\T_{\varphi}$ is a homeomorphism (this relies on the notion
of {\em recognizability} or invertibilty for such substitutions--see \cite{mosse}
and
\cite{sol}). Each
$T \in \T_{\varphi}$  is uniformly recurrent under the flow and has dense orbit
(i.e.,   the flow is minimal on $\T_{\varphi}$).  It follows that 
$\T_{\varphi}$ is a continuum.  

Recall that a {\em composant} of a point $x$ in a topological space $X$ is the
union of the proper compact  connected subsets of $X$ containing $x$. If $\vp$ is
a primitive substitution,  
composants and arc components in $\T_\vp$ are identical; in this case we use
the terms interchangeably.  For any  substitution
$\varphi$, the arc components of the tiling space
${\mathcal T}_\varphi$ coincide with the orbits of the natural flow (translation) on
${\mathcal T}_\varphi$.  In particular, if ${\mathcal C}$ is an arc component 
 of ${\mathcal T}_\varphi$, $$\{ \underline{w}(T): T \in {\mathcal C}\} =
[\underline{w}(T)],$$ where $[
\underline{w}(T)]$ is the shift class of 
$\underline{w}(T)$.  We also call  $[
\underline{w}(T)]$ the {\em pattern of the arc component (or composant) of
$T$}. 

Tilings $T, T' \in \T_\vp$ are {\em forward asymptotic} if $\lim_{t \arrow
\infty} dist(T - t, T' - t) = 0$. Equivalently, $T =\{T_i\}_{i=-\infty}^{\infty},
T' = \{T'_i\}_{i=-\infty}^{\infty}  $ are {\em forward asymptotic} if there are
$N, M
\in
\Z$ so that $T_{N+k} = T'_{M+k} $ for all $k \geq 0$. 
Composants are {\em forward asymptotic} if they contain forward asymptotic
tilings.  Backward asymptotic tilings and composants are defined similarly.

 Given a primitive substitution  $\varphi : \A \arrow \A^*$ with $card(\A) =
d \geq 2$, and left Perron-Frobenius eigenvector $ \omega_L = (\omega_1, \ldots,
\omega_d)$, define
$R_\vp =
\vee^d_{i = 1} S_i$ as a wedge of $d$ oriented circles $S_1, \ldots, S_d$
with the  circumference of $S_i = \omega_i$, and let $f_\varphi : R_\vp
\arrow R_\vp$ be the `linear' map, with expansion constant $\lambda$,
that follows the pattern
$\varphi$.  That is, if $\varphi(i) = {i_1}{i_2}
\ldots {i_{k(i)}}$,   then $f_\varphi $ maps the circle $S_i$ around the circles
$S_{i_1},  \ldots, S_{i_{k(i)}}$, in that order, preserving orientation
and stretching distances locally by a factor of $\lambda$.   We call
$f_\varphi$  {\em the map of the rose associated with
$\varphi$}. 

As with the case for tiling spaces, there are situations
when we wish to define $R_\vp$ but $\vp$ is not
primitive.  In these situations, if $A_\vp$ does not have  a unique left
eigenvector,  any positive left eigenvector will suffice to define the
circumferences of the circles. 

If $f: X \to X$ is a map of a compact connected metric space $X$, then the {\em
inverse limit space} with single bonding map $f$
 is the space $$\inv f  = \{(x_0,x_1,\ldots) :
f(x_i) = x_{i-1}\;{\rm for}\; i=1,2,\ldots\}$$ with metric
$$\underline{d}(\underline{x}, \underline{y}) = \Sigma_{i \geq 0}\frac{d(x_i,
y_i)}{2^i};$$ 
$\hat{f}:\inv f\to \inv f$ will denote the natural (shift)
homeomorphism $$\hat{f}(x_0,x_1,\ldots)=(f(x_0),x_0,x_1,\ldots).$$

We now define a second model of the tiling space called {\em strand space}.  
A {\em strand}, $\gamma$, in $\R^d$ is a collection of {\em segments},
(sometimes called {\em edges}), $\gamma = \{S_n\}_{n = N}^M$, with each
segment $S_n$ a translate of a unit interval parallel to a coordinate axis: 
if $\{e_1, \ldots, e_d\}$ is the
standard basis for $\R^d$, and  $I_i := \{te_i: 0
\leq t \leq 1\}$, then $S_n = I_{i_n} + v_n$ for some $i_n \in \{1, \ldots,
d\}$ and $v_n \in \R^d$.  Moreover, we require that $\gamma$ be `connected'
in the sense that $v_{n+1} = v_n + e_{i_n}$, $i = N, \ldots, M-1$.  Two
strands are equal if they are identical as collections of segments. In
particular, if
$S_n = S'_{n + l}$ for $n = N,
\ldots, M$, then $\{S_n\}_{n = N}^M = \gamma = \gamma' = \{S'_n\}_{n = N +
l}^{M+l}$.  A strand $\gamma = \{S_n\}_{n = -\infty}^\infty$ is said to be 
{\em bi-infinite}.  Let $$\F^d = \{ \gamma: \gamma \textrm{ is a
bi-infinite strand in } \R^d\}.$$

The endpoints of segments in a strand are called the {\em vertices}, and if $S_n =
I_{i_n} + v_n$, then $\min S_n := v_n$, $\max S_n := v_n + e_{i_n}$. If
$\gamma = \{S_n\}_{n = N}^M$ and $\gamma' = \{S'_n\}_{n = N'}^{M'}$ are two
strands with $M, N' < \infty$ and $\max S_M = \min S'_{N'}$, we can
concatenate $\gamma$ and $\gamma'$ (that is, union and reindex) to obtain a
single longer strand $\gamma \cup \gamma'$.

Given a substitution $\vp$ on the alphabet $\A = \{1, \ldots, d\}$, define
$\Phi(\{I_i\})$ to be the strand $\{S_{n,i}\}_{n = 1}^{k(i)}$ with $S_{n,i}
= I_{i_n} + (\Sigma_{j=1}^{n-1} e_{i_j})$, where $\vp(i) = i_1 \ldots
i_{k(i)}$.  That is, $\Phi$ applied to the  singleton strand $I_i$ is the
strand with the origin as initial vertex that `follows the pattern' of the
word $\vp(i)$.  Now extend $\Phi$ to arbitrary singleton strands by $$
\Phi(\{I_i + v\}) := \{I_{i_n} + (\Sigma_{j=1}^{n-1} e_{i_j}) + Av\}_{n =
1}^{k(i)}  
$$ and to arbitrary strands by concatenation:  $$\Phi(\{S_n\}_{n = N}^M) :=
\cup_{n = N}^M \Phi( \{S_n\}).$$  For $R \in \R$, let $\F^d_R$ denote the
subspace
 of $\F^d$ consisting of those strands that lie in a cylinder of
radius $R$ centered on $E^u$. If $\vp$ is weak Pisot, for sufficiently
large
$R$, 
$\F^d_R$ is mapped into itself by $\Phi$: choose $R_0$  so that
$\Phi(\F^d_{R_0}) \subset \F^d_{R_0}$.  Define $$\F^d_{\vp} =
\{\gamma = \{S_n\}_{n = -\infty}^\infty \in \F^d_{R_0}: \textrm{ if } S_k
\cap E^s \neq \emptyset, \textrm{ then } i_{k-1}i_ki_{k+1} \in \cL_\vp\}.
$$
The {\em strand space of $\vp$} is $$\T_\vp^S := \cap_{n \geq 0}
\Phi^n(\F_\vp^d).$$

A metric can be defined on  $\T^S_\vp$ that has the property:  the distance
between $\gamma = \{S_n\}_{n = -\infty}^\infty$ and $\gamma' = \{S'_n\}_{n =
-\infty}^\infty$ is small if there is $v \in \R^d$, $|v|$ small, and $N \in
\N$, $N $ large, so that $S_n = S'_n + v$ for $n = -N, \ldots, N$ (where the
indexing is such that $S_0 \cap E^s \neq \emptyset $).

It is proved in \cite{mjarek} that if $\vp$ is a strong Pisot substitution,
then $\Phi:  \T^S_\vp \arrow \T^S_\vp $ is a homeomorphism
(referred to as the $\Z$-action), $(\gamma, t) \mapsto \gamma + t \omega_R$
defines a flow on $\T^S_\vp$ (referred to as the $\R$-action), and there is
a homeomorphism $h : \T^S_\vp \arrow \T_\vp$ that conjugates to $\Z$- and
$\R$-actions on these spaces ($h$ is just the projection of a strand onto
$E^u \simeq \R$ along $E^s$).

In what follows, we refer to {\em rewriting} a substitution.  When this
term is used without qualification, we are referring to rewriting as developed in
\cite{dur}, to which we refer the reader.  We discuss the construction in some
detail and provide examples in
\cite{bd1}.  We also include an example in \S \ref{sec-prox} of this paper.

Finally, $X \simeq Y$ will mean that $X$ and $Y$ are homeomorphic.

\section{The No Cycle Condition}\label{sec-cohomol}
The main result of this section is Corollary \ref{cohomcor}, from which we
eventually deduce that if
$\vp$ is strong Pisot, and  $\tv$ is  a {\em one-cut rewriting of
$\vp$} (see \S \ref{sec-prox}), then $\tv$ is weak Pisot. A
reader willing to accept Corollary \ref{cohomcor}  can safely proceed to \S
\ref{sec-prox}. 

Given a substitution  $\vp$, 
  let
$\vp^+$, $\vp^-: \A \arrow \A$ be defined by:  $$\vp^+(a) = b \textrm{ if }
\vp(a) = b  \ldots,$$  $$\vp^-(a) = c \textrm{ if } \vp(a) =
 \ldots  c.$$ Let $S^+ = \cap_{n \geq 0} (\vp^+)^n(\A)$ and $S^- = \cap_{n
\geq 0} (\vp^-)^n(\A)$ be the eventual ranges of $\vp^+$ and $\vp^-$, and let
$\cP_\vp = \{(a, b) \in S^- \times S^+: ab \in \cL_\vp\}$.  Define an equivalence
relation
$\sim$ on $\cP_\vp$ by $(a, b) \sim (c, d)$ if $a = c$ or $b = d$ and extending
by transitivity.  A {\em cycle } in $\cP_\vp$ consists of a string of
equivalences $(a_1, a_2) \sim (a_3, a_2) \sim (a_3, a_4)
\sim \ldots \sim (a_1, a_{2n})$ with $a_1 \neq a_3 \neq \ldots \neq
a_{2n-1}$, $a_2 \neq a_4 \neq \ldots
\neq a_{2n}$, and $n \geq 2$.

\

{\parindent=0pt {\bf  No Cycle Condition:}}
The substitution $\vp$ has   {\em no cycles of periodic
words} if $\cP_\vp$ has
no cycles.

\

Note that if $(a, b)$, $(c, b) \in \cP_\vp$, then
the bi-infinite words obtained by iterating on $a.b$ and $a.c$ ($\lim_{n
\arrow \infty} \vp^n(a) . \vp^n(b)$ and $\lim_{n
\arrow \infty} \vp^n(c) . \vp^n(b)$ respectively) are periodic under $\vp$ and
represent asymptotic composants.  That is, a cycle in $\cP_\vp$ represents a
cycle of asymptotic composants $\{\C_1, \ldots, \C_n\}$ ($\C_1,
\C_2$ are forward asymptotic, $\C_2, \C_3$ are
backward asymptotic, etc.) of a particular sort.

The existence of a cycle of asymptotic composants is a topological property
of a space.  However, whether a cycle of asymptotic composants in $\T_\vp$ is
associated with a cycle in $\cP_\vp$ is combinatorial and an 
artifact of the symbolic presentation of the tiling space.  For instance, the
Fibonacci substitution  ($\psi(1) = 12, \psi(2) = 1$) has a cycle of two
composants asymptotic in both directions but 
$\cP_\psi$ has no cycles.  On the other hand, the  Morse-Thue substitution
$\vp$  defined  by 
$\vp(1) = 12, \vp(2) =  21$ has a cycle of four asymptotic composants, all
associated with the four bi-infinite words periodic under $\vp$: $\ldots 1.1
\ldots$, $\ldots 1.2
\ldots$, $\ldots 2.1 \ldots$, $\ldots 2.2 \ldots$, which determine a cycle in
$\cP_\vp$.  There are several ways to rewrite
$\vp$  to obtain a proper substitution $\vp'$ for which the tiling spaces $\T_\vp$
and
$\T_{\vp'}$ are homeomorphic; since $\vp'$ has a single bi-infinite periodic
word, its cycle of asymptotic composants is no longer associated with a cycle
in $\cP_{\vp'}$.   

The cycle of asymptotic composants in the tiling
space $\T_\vp$ (where $\vp$ is the Morse-Thue substitution)  is associated with a
nontrivial element of the cohomology of $\T_\vp$ that is periodic under
the action induced by inflation and substitution.  That is, inflation and
substitution has a root of unity eigenvalue on the level of cohomology.  Thus any
proper rewriting of $\vp$ has an abelianization that has a root of unity
eigenvalue and hence is non-Pisot.  The non-Pisot nature of the rewriting
destroys the connection between the symbolic and the geometrical aspects of
proximality that we exploit to obtain the results of this paper. In the case of
the Fibonacci substitution, there is also a cohomology element associated with the
cycle of asymptotic composants, but it is trivial, and does not correspond to a
root of unity eigenvalue--in this case, our program can proceed.

 For these reasons, in this paper we require that the substitutions we consider
have no cycles of periodic words.
 
  Let
$\cP_\vp = \cP$ be defined as above.
\begin{lemma} \label{singleeq} If $\vp$ is weak Pisot, then $\cP$ consists of
a single equivalence class.
\end{lemma}

Proof:  Suppose that there are $k$ equivalence classes in $\cP$, where $k
\geq 2$. Using the elements of $\cP$ to define starting and stopping rules,
we obtain a rewriting $\tilde{\vp}$ of $\vp$ with alphabet $\tilde{\A}$
consisting of certain elements of $\cL(\vp)$.  Let $v_1, \ldots, v_k$ denote
the equivalence classes of
$\sim$ in $\cP$.  Define a graph $G_{\tilde{\vp}} = G$ with vertices  $v_1,
\ldots, v_k$  and edges labeled by elements of $\tilde{\A}$:  the edge
labeled $b \ldots a
\in \tA$ starts at vertex $v_i = [(\hspace{.1in}, b)]$ and ends at vertex
$v_j = [(a, \hspace{.1in})]$.  Let
$f_{\tv} : G \arrow G$ follow pattern $\tv$, $R_\vp$ be the rose with
vertex
$v$ and edges labeled by $\A$, and 
$f_\vp : R_\vp \arrow R_\vp$ follow pattern $\vp$. Let $\rho : \tA \arrow
\A^*$ be the natural morphism, and  $f_\rho: G_{\tv} \arrow R_\vp$ be the
map following pattern $\rho$.  The maps $f_\vp$, $f_{\tv}$ and $f_\rho$ satisfy
$f_\rho
\circ f_{\tv} = f_\vp
\circ f_\rho$.  There is then an induced surjection $\hat{f}_\rho: \inv
f_{\tv}
\arrow
\inv f_\vp$.  It is easily checked that
$\hf_\rho$ is one-to-one everywhere except at a single point:  if
$\underline{v}^i := (v_i, v_i, \ldots)$ for $i = 1, \ldots, k$ and
$\underline{v} := (v, v, \ldots)$, then 
$(\hf_\rho)^{-1}(\underline{v}) = \{\underline{v}^1, \underline{v}^2
, \ldots, \underline{v}^k \} := \V. $  Thus
$\inv f_\vp \simeq (\inv f_{\tv})/ \V$.

We have the commuting diagram derived from the long exact sequences of pairs
(the zero-th level is reduced and the coefficients are $\Q$):

{\parindent=0pt
\begin{picture}(20, 140)(-115, -110)

\put(-23, 20){$0 \arrow 0 \arrow \check{H}^1(\inv f_\vp, \{\underline{v}\})
\arrow \check{H}^1(\inv f_\vp ) \arrow 0$}

\put(-50,-20){$0 \arrow 0  \arrow \check{H}^1(\inv f_\vp, \{\underline{v}\})
\arrow \check{H}^1(\inv f_\vp ) \arrow 0$}

\put(-39, -60){$0 \arrow \check{H}^0(\V)$}
\put(21, -57){$\line(1,0){5}\  \vector(1,0){8} $} 
\put(40, -60){$\check{H}^1(\inv f_{\tilde{\vp}}, \V)$}
\put(115, -57){$\line(1,0){5}\  \vector(1,0){8} $} 
\put(137, -60){$ \check{H}^1(\inv
f_{\tilde{\vp}}) \arrow 0$}

\put(-70,-100){$0 \arrow \check{H}^0(\V)  \arrow
\check{H}^1(\inv f_{\tilde{\vp}}, \V)
\arrow \check{H}^1(\inv f_{\tilde{\vp}} ) \arrow 0$}

\put(60, 12){\vector(-1,-1){20}}
\put(33, 3){$\hat{f}^*_\vp$}

\put(150, 12){\vector(-1,-1){20}}
\put(123, 3){$\hat{f}^*_\vp$}

\put(150, -68){\vector(-1,-1){20}}
\put(145,-78){$\hat{f}^*_{\tilde{\vp}}$}
\put(60,-68){\vector(-1,-1){20}}
\put(58,-78){$\hat{f}^*_{\tilde{\vp}}$}
\put(-0,-66){\vector(-1,-1){20}}
\put(-28,-76){$\hat{f}^*_{\tilde{\vp}}$}

\put(5,15){\line(0,-1){20}}
\put(5,-24){\vector(0,-1){22}}

\put(75,12){\line(0,-1){20}}
\put(75,-25){\vector(0,-1){22}}
\put(78, 0){$f^*_\rho$}
\put(160,12){\line(0,-1){22}}
\put(160,-22){\vector(0,-1){25}}
\put(162, 0){$f^*_\rho$}

\put(28,-28){\vector(0,-1){60}}
\put(122,-28){\vector(0,-1){60}}

\end{picture}

The} bottom rows split since all the homomorphisms are  linear maps of
finite dimensional vector spaces.  Let $   \hat{f}_{\tilde{\vp} }|_\V
:= 
\hat{f}_{\tilde{\vp}, 1 }$. Then $\hat{f}^n_{\tilde{\vp}, 1 } = id$ for $n = k!$ 
We extract the commuting square

\begin{picture}(25, 100)(-120, -55)

\put(-24, -7){$h$}
\put(115, -7){$h$}
\put(44, 25){$\hf_\vp^*$}
\put(30, -28){$\hat{f}^*_{\tilde{\vp}, 1 } \oplus \hat{f}^*_{\tilde{\vp}, 2
}$}

\put(-55, 17){$\check{H}^1(\inv f_\vp, \{\underline{v}\})$}
\put(85, 17){$\check{H}^1(\inv f_\vp, \{\underline{v}\})$}
\put(-75, -38){$\check{H}^0(\V) \oplus \check{H}^1(\inv f_{\tv})$}
\put(85, -38){$\check{H}^0(\V) \oplus \check{H}^1(\inv f_{\tv})$}

\put(-10,12){\vector(0, -1){37}}
\put(34,-35){\vector(1,0){45}}
\put(130,12){\vector(0,-1){37}}
\put(34,20){\vector(1,0){45}}

\end{picture}
\\of vector space isomorphisms.  Since $\dim \check{H}^0(\V) = k-1 \geq 1   $
and all eigenvalues of $\hat{f}^*_{\vp, 1 } $ are roots of unity, 
$\hat{f}^*_{\vp} $ has roots of unity eigenvalues.  

Finally, by  continuity of the \v{C}ech theory,
$$\check{H}^1(\inv f_\vp,
\{\underline{v}\})
\\ \simeq
\check{H}^1(\inv f_\vp)
 \simeq \dir (A^{tr}_\vp: \Q^d \arrow \Q^d) \simeq
ER_\vp,$$ 
where
$ER_\vp$ is the eventual range of $A^{tr}_\vp$, and $\hf_\vp^*$ is conjugated
to
$A^{tr}_\vp:ER_\vp \arrow ER_\vp$ by this isomorphism. Thus
$A_\vp$ also has root of unity eigenvalues, contradicting the assumption that
$\vp$ is weak Pisot. 
\qed

\

Remark: The proof of the above lemma only requires that  $\vp$ is
 hyperbolic.

Suppose that $\psi_t: X \arrow
X$ is a flow on the compact metric space
$X$, and that  for each $j = 1, \ldots, m$, $X_{1, j}, \ldots, X_{n_j, j}$ 
are asymptotic (forward, say) under $\psi_t$: that
is,
$d(\psi_t(X_{i,j}), \psi_t(X_{l,j})) \arrow 0$ as $t \arrow \infty$ for $i, l
\in
\{1, \ldots, n_j\}$ and $j \in \{i, \ldots, m\}$. Define $\overline{X}$ to be
the quotient
$\overline{X} = X/\sim$, where $\psi_t(X_{i,j}) \sim \psi_t(X_{l,j}) $ for $t
\geq 0$, $i, l
\in
\{1, \ldots, n_j\}$ and $j \in \{i, \ldots, m\}$.

\begin{lemma} \label{inducemap} The quotient map $p : X \arrow \overline{X}$
induces an isomorphism $p^*: \check{H}^k(\overline{X}) \arrow
\check{H}^k(X)$ for all
$k$.
\end{lemma}

Proof: To avoid excessive notation, we assume $m = 1$. Let
$\overline{\psi}_t$ be the semi-flow on $\overline{X}$ defined by
$\overline{\psi}_t(p(x)) = p(\psi_t(x))$ for $t \geq 0$.  Define $$\inv
\overline{\psi}_t := \{ \gamma: \R \arrow \overline{X} :
\overline{\psi}_t(\gamma(s)) = \gamma(t + s) \textrm{ for all } t \geq 0, s
\in \R\}$$ with the compact-open topology.  Because the glued arcs are
asymptotic,
$\overline{X}$ is compact and Hausdorff, hence $ \inv
\overline{\psi}_t$ is compact and Hausdorff.  Let $\hat{p}: X \arrow \inv
\overline{\psi}_t$ be defined by $\hat{p}(x) =
p(\psi\underline{\hspace{.07in}}(x)) :
\R
\arrow X$.  Then $\hat{p}$ is clearly continuous.  

To see that $\hat{p}$ is a surjection, consider $\psi_t(
p^{-1}(\gamma(-t)))$ for $\gamma \in \inv
\overline{\psi}_t$.  Either this is a singleton or $\gamma(-t) \in \alpha
:= \{p(\psi_s(x_i)): s \geq 0\}$.  In
the latter case,
$p^{-1}(\gamma(-t)) = \{\psi_s(x_i) : i = 1, \ldots, n\}$ for some $s
\geq 0$.  Then $\psi_t(p^{-1}(\gamma(-t))) = \{\psi_{t+s}(x_i): i = 1,
\ldots, n\}$ for some $s \geq 0$, and the diameter of this set goes to 0 as
$t \arrow \infty$, uniformly in $s \geq 0$.  In any event,
$diam ( \psi_t
(p^{-1}(\gamma(-t)))) \arrow 0$ as $t \arrow \infty$.  Note also that
$$p(\psi_t (p^{-1}(\gamma(-t)))) = \overline{\psi}_t(\gamma(-t)) =
\gamma(0),$$ so $x := \lim_{t \arrow \infty}\psi_t
(p^{-1}(\gamma(-t)))$ is well-defined.  Moreover, $\hat{p}(x) =
p(\psi\underline{\hspace{.07in}}(x))$ and
\begin{eqnarray*} 
p(\psi_t(x))& =  p \psi_t(\lim_{s \arrow \infty}
\psi_s(p^{-1}(\gamma(-s)))) \\ & =  \lim_{s \arrow \infty}
p (\psi_{t+s}(p^{-1}(\gamma(-s)))) \\&= \lim_{s \arrow \infty}
\overline{\psi}_{t+s}(\gamma(-s)) \hspace{.6in}  \\ &= \lim_{s
\arrow
\infty}
\gamma(t) \hspace{1.15in} \\&= \gamma(t). \hspace{1.64in}
\end{eqnarray*}  
Thus $\hat{p}$ is surjective.

If $\hat{p}(x) = \hat{p}(y) = \gamma$, then $p(\psi_t(x)) = p(\psi_t(y)) =
\gamma(t)$ for all $t$.  Then $\psi_t(x), \psi_t(y) \in
p^{-1}(\gamma(t))$ and $x, y \in \psi_{-t}(p^{-1}(\gamma(t)))$ for
all $t$.  But $diam (\psi_{-t}(p^{-1}(\gamma(t)))) \arrow 0$ as $t \arrow
- \infty$, so $x = y$.  Since $X$ is compact and $\inv \overline{\psi}_t$ is
Hausdorff, $\hat{p}$ is a homeomorphism.

Now let $\overline{\psi}_1: \overline{X} \arrow \overline{X}$ be the time
one map of $\overline{\psi}$.  If $h : \inv \overline{\psi}_1 \arrow
\inv \overline{\psi}_t$ is defined by $h(\overline{x}_0, \overline{x}_1,
\ldots) =
\gamma$, where $\gamma(t) = \overline{\psi}_{n+t}(\overline{x}_n)$ for $n
\geq -t$, then
$h$ is a homeomorphism.  By  continuity of the \v{C}ech theory,
$$\check{H}^k(X) \stackrel{\hat{p}^*}{\simeq}  \check{H}^k(\inv
\overline{\psi}_t) \stackrel{h^*}{\simeq} \check{H}^k(\inv
\overline{\psi}_1) \simeq \dir \overline{\psi}_1^*:
\check{H}_1(\overline{X}) \arrow \check{H}_1(\overline{X}).$$
Finally, since $\overline{\psi}_1$ is homotopic to 
$
\overline{\psi}_0 = id$ (where $\overline{\psi}_t$, $0 \leq t \leq 1$,
provides the homotopy), $\overline{\psi}_1^* = id$ and $\dir
\overline{\psi}_1^* \simeq \check{H}_1(\overline{X})$ (by projection onto the
first coordinate). The composition of these isomorphisms is
$p^*$.
\qed

 \ 

Recall that $\Phi: \T_\vp \arrow \T_\vp$ is the inflation and
substitution map.

\begin{prop} Suppose that $\vp$ is weak Pisot and has no cycles of periodic
words.  Then the linear map on $\dir (A^{tr}_\vp: \Q^d \arrow \Q^d)$ induced by
$A^{tr}_\vp$ is conjugate to the isomorphism $F^*_\vp: \check{H}^1 (\T_\vp) \arrow
\check{H}^1 (\T_\vp) $.
\end{prop}

Proof:  By Lemma \ref{singleeq}, $\cP$ consists of a single equivalence
class.  For each $(a, b)
\in \cP$, let $T_{(a, b)}$ be the corresponding (periodic under $\Phi$)
tiling in
$\T_\vp$.  Let 
\begin{eqnarray*} \cP^+ &:= \{ b \in \A: \textrm{ there are
} a
\neq c
\textrm{ with } (a, b), (c, b) \in \cP\} \\ &= \{b_1, \ldots,
b_k\}. \hspace{2.45in}
\end{eqnarray*}
  For each $b_j \in \cP^+$, let $\{a_{i,j}\}_{i = 1}^{l_j}$ be a list of
the letters for which $(a_{i,j}, b_j) \in \cP$.  Let
$\overline{X}^+$ be the quotient space obtained from $\T_\vp$ by the
identifications
$T_{(a_{i,j}, b_j) } - t \sim T_{(a_{m,j},
b_j) } - t$ for $t \geq 1$, $1 \leq i, m \leq l_j$, and $j = 1, \ldots, k$. 
By Lemma \ref{inducemap},   the quotient map $p_+: \T_\vp \arrow
\overline{X}^+$  induces a vector space  isomorphism 
$p_+^*:
\check{H}^1(\overline{X}^+) \arrow \check{H}^1(\T_\vp)$ (here and in what
follows, coefficients are $\Q$).   Furthermore, the inflation and
substitution map,
$\Phi$, induces a map
$F_+: \overline{X}^+ \arrow \overline{X}^+$ so that 

{\parindent=0pt
\begin{picture}(25, 100)(-140, -55)

\put(0, -7){$p_+^*$}
\put(95, -7){$p_+^*$}
\put(44, 25){$F_+^*$}
\put(44, -28){$\Phi^*$}

\put(-10, 17){$\check{H}^1(\overline{X}^+)$}
\put(85, 17){$\check{H}^1(\overline{X}^+)$}
\put(-10, -38){$\check{H}^1(\T_\vp) $}
\put(85, -38){$\check{H}^1(\T_\vp)$}

\put(15,12){\vector(0, -1){37}}
\put(34,-35){\vector(1,0){45}}
\put(110,12){\vector(0,-1){37}}
\put(34,20){\vector(1,0){45}}

\end{picture}

is a} commuting diagram of vector space isomorphisms.

Let $S: \overline{X}^+ \arrow [0, \infty)$ be a continuous map with
$$S^{-1}(0) = \{ [T_{a_{i, j},b_j} - 1]: j = i, \ldots k\}.$$  That is,
$S$ vanishes exactly at the branch points of $\overline{X}^+$.  Moreover,
using the local product structure of $\T_\vp$, we may choose $S$ so
that if
$(a, b), (a, c) \in \cP$, then $S(p_+(T_{(a, b)} + t)) = S(p_+(T_{(a, c)} +
t)) $ for all $t \geq 0$.  Let $\tau: \T_\vp \times \R \arrow \R$ be the
solution to 
$$ \frac{d \tau} {dt} (T, t) = S \circ p_+ (T - \tau(T, t)), \hspace{.5in}
\tau(T, 0) = 0.$$  
Then $(T, t) \mapsto T - \tau(T, t)$ is a flow on
$\T_\vp$ that descends to a flow $\psi_t$ on $\overline{X}^+ $ with rest
points at exactly the branch points $\{ [T_{a_{i, j},b_j} ]: j = i, \ldots
k\}$. 

Now let \begin{eqnarray*} \cP^- &:= \{a \in \A: \textrm{ there are } b \neq c
\in
\A
\textrm{ with } (a, b), (a, c) \in \cP \}\\ &= \{a_1, \ldots,
a_l\}. \hspace{2.8in}
\end{eqnarray*}
  For each
$a_i
\in \cP^-$, let  $\{b_{i,j}\}_{j = 1}^{k_i}$ be a
list of the letters for which $(a_{i}, b_{i, j}) \in \cP$. Let
$\overline{X}$ be the quotient space obtained from $\overline{X}^+$ by the
identifications 
$\psi_{-t}(p_+(T_{(a_i,
b_{i,j}) })) \sim \psi_{-t}(p_+(T_{(a_{i},
b_{i,m}) } ))$ for $t \geq 1$, $1 \leq j, m \leq k_i$, and $i = 1, \ldots,
l$, with quotient map  $p_-: 
\overline{X}^+ \arrow \overline{X}$. There is an
induced map $F: \overline{X} \arrow \overline{X}$ with $p_- \circ F_+ = F
\circ p_-$ and,  from Lemma
\ref{inducemap},
we have a commuting diagram of vector space isomorphisms:

{\parindent=0pt
\begin{picture}(25, 100)(-140, -55)

\put(0, -7){$p_-^*$}
\put(95, -7){$p_-^*$}
\put(44, 25){$F^*$}
\put(44, -28){$F_+^*$}

\put(-10, 17){$\check{H}^1(\overline{X})$}
\put(85, 17){$\check{H}^1(\overline{X})$}
\put(-10, -38){$\check{H}^1(\overline{X}^+) $}
\put(85, -38){$\check{H}^1(\overline{X}^+)$}

\put(15,12){\vector(0, -1){37}}
\put(34,-35){\vector(1,0){45}}
\put(110,12){\vector(0,-1){37}}
\put(34,20){\vector(1,0){45}}

\end{picture}

 Letting} $p = p_- \circ p_+$, and combining the above diagrams, we have

{\parindent=0pt
\begin{picture}(25, 100)(-140, -55)

\put(0, -7){$p^*$}
\put(95, -7){$p^*$}
\put(44, 25){$F^*$}
\put(44, -28){$\Phi^*$}

\put(-10, 17){$\check{H}^1(\overline{X})$}
\put(85, 17){$\check{H}^1(\overline{X})$}
\put(-10, -38){$\check{H}^1(\T_\vp) $}
\put(85, -38){$\check{H}^1(\T_\vp)$}

\put(15,12){\vector(0, -1){37}}
\put(34,-35){\vector(1,0){45}}
\put(110,12){\vector(0,-1){37}}
\put(34,20){\vector(1,0){45}}

\put( 160, 0){(1)}

\end{picture}
 
Let}  $q: \T_\vp \arrow R_\vp$ denote the map that takes a tiling to
the location of its origin (that is, if the $i^{th}$ petal $P_i$ of
the rose $R_\vp$ is identified with the $i^{th}$ prototile $[0, \lambda_i]$,
and
$[0,
\lambda_i] - t$ is the tile of $T$ containing 0, then $q(T) = t \in P_i$).
Let $f_\vp: R_\vp \arrow R_\vp$ be
the rose map and
$\hat{q}: \T_\vp \arrow \inv f_\vp$  the continous surjection defined by
$$\hat{q}(T) = (q(T), q(f^{-1}_\vp(T), q(f^{-2}_\vp(T), \ldots).$$  There is then
$\overline{q}:
\overline{X} \arrow
\inv f_\vp$ so that $\overline{q} \circ p = \hat{q}$ and  $\overline{q}
\circ F = \hat{f}_\vp \circ \overline{q}$.  Let 
\begin{eqnarray*}
B &=& \{T_{(a_{i,j},
b_j) } - t : 0 \leq t \leq 1;  i = 1, \ldots, l_j; j = 1,
\ldots, k\} \\ & &\  \cup \  \{T_{(a_{i},
b_{i, j}) } + t : 0 \leq t \leq 1;  j = 1, \ldots, k_i; i = 1,
\ldots, l\}
\end{eqnarray*}
 (where $S$ must be adjusted so that $S \circ p_+(T_{(a_{i},
b_{i, j}) } + t) := 1$ for  $0 \leq t \leq 1;  j = 1, \ldots, k_i; i = 1,
\ldots, l$), let $\overline{B} = p(B)$, and let $\hat{B} =
\overline{q}(\overline{B}) = \hat{q}(B)$.
The fact that $\cP$ is one equivalence class means that $\overline{B}$
is connected, and the assumption that $\cP$ has no cycles means that
$\overline{B}$ is acyclic; that is, $\overline{B}$ is a finite tree and is
contractible.  Since $F(\overline{B}) \supset \overline{B}$,  we may
homotop $F$ to $G$ with $G(\overline{B}) = \overline{B}$.  Similarly, we
may homotop  
$\hat{f}_\vp$ to $g$ so that  $g(\hat{B})= \hat{B}$ and $\overline{q} \circ
G = g \circ  \overline{q}$.  One can further check that $\overline{q}$ is
one-to-one except on $\overline{B}$.  We have the commuting diagram:

{\parindent=0pt
\begin{picture}(25, 100)(-140, -55)

\put(-17, -7){$\overline{q}^*$}
\put(115, -7){$\overline{q}^*$}
\put(44, 25){$g^*$}
\put(44, -28){$G^*$}

\put(-40, 17){$\check{H}^1(\inv f_\vp, \hat{B})$}
\put(85, 17){$\check{H}^1(\inv f_\vp, \hat{B})$}
\put(-30, -38){$\check{H}^1(\overline{X}, \overline{B})$}
\put(95, -38){$\check{H}^1(\overline{X}, \overline{B})$}

\put(-5,12){\vector(0, -1){37}}
\put(34,-35){\vector(1,0){45}}
\put(127,12){\vector(0,-1){37}}
\put(34,20){\vector(1,0){45}}

\end{picture}

 Since}  $\overline{q}: \overline{X}/\overline{B} \arrow \inv f_\vp/ \hat{B}$
is a homeomorphism, $\overline{q}^*$ is an isomorphism.  Thus 

{\parindent=0pt
\begin{picture}(25, 100)(-140, -55)

\put(-5, -7){$\overline{q}^*$}
\put(104, -7){$\overline{q}^*$}
\put(44, 25){$g^*$}
\put(44, -28){$G^*$}

\put(-25, 17){$\check{H}^1(\inv f_\vp)$}
\put(85, 17){$\check{H}^1(\inv f_\vp)$}
\put(-10, -38){$\check{H}^1(\overline{X})$}
\put(95, -38){$\check{H}^1(\overline{X})$}

\put(7,12){\vector(0, -1){37}}
\put(34,-35){\vector(1,0){45}}
\put(118,12){\vector(0,-1){37}}
\put(34,20){\vector(1,0){45}}

\end{picture}

is} a commuting diagram of isomorphisms. Replacing $g^*$ by $\hat{f}_\vp^*$ and
$G^*$ by $F^*$ (to which they are equal, respectively), and combining this
last with  diagram (1), we have the commuting diagram of vector space
isomorphisms:

\begin{picture}(45, 100)(-125, -55)

\put(-3, -7){$q^*$}
\put(95, -7){$q^*$}
\put(44, 25){$\hat{f}_\vp^*$}
\put(44, -28){$\Phi^*$}

\put(-25, 17){$\check{H}^1(\inv f_\vp)$}
\put(85, 17){$\check{H}^1(\inv f_\vp)$}
\put(-10, -38){$\check{H}^1(\T_\vp) $}
\put(85, -38){$\check{H}^1(\T_\vp)$}

\put(10,12){\vector(0, -1){37}}
\put(34,-35){\vector(1,0){45}}
\put(110,12){\vector(0,-1){37}}
\put(34,20){\vector(1,0){45}}

\end{picture}

By continuity of the \v Cech theory, there is an isomorphism between
$\check{H}^1(\inv f_\vp)$ and $\dir f_\vp^*$  (where $ f_\vp^*:
H^1(R_\vp) \arrow H^1(R_\vp)$) that conjugates $\hat{f}_\vp^*$ with 
$\widehat{f_\vp^*} : \dir f_\vp^* \arrow \dir f_\vp^*$ defined by
$$\widehat{f_\vp^*} ([(\gamma, n]) = [(f_\vp^*(\gamma), n)].$$ 
Finally, identifying $H^1(R_\vp)$ with the dual of $H_1(R_\vp)$ and
choosing as ordered basis for  $H_1(R_\vp)$ the oriented petals of
$R_\vp$, the matrix for $(f_\vp)_*$ is $A_\vp$, and $A^{tr}_\vp$ represents
$f_\vp^*$ (where $A^{tr}$ is the transpose of $A$).  Thus $\widehat{f_\vp^*}$ is
conjugate with
$A^{tr}_\vp : 
\dir (A^{tr}_\vp: \Q^d \arrow \Q^d) $.
\qed

\begin{corollary} \label{cohomcor} If $\vp$ is weak Pisot and has no
 cycles of periodic words, and $\psi$ is a proper substitution with the property
that
$\Phi$ and
$\Psi$ are conjugate, then
$\psi$ is weak Pisot.

\end{corollary}

Proof:  Since $\psi$ is proper, there is a homeomorphism
of $\T_{\psi}$ with
$\inv f_{\psi}$ that conjugates $F_{\psi}$ with $\hat{f}_{\psi}$.  We
have a conjugacy between 
$A^{tr}_{\vp}$ on $\dir (A^{tr}_\vp: \Q^d \arrow \Q^d)$ and $A^{tr}_{\psi} $
on $\dir (A^{tr}_\psi: \Q^d \arrow \Q^d)$.  But  $\dir (A^{tr}_\vp: \Q^d
\arrow \Q^d)$ and
$\dir (A^{tr}_\psi: \Q^d \arrow \Q^d)$ are naturally isomorphic with the
eventual ranges of
$A^{tr}_\vp : \Q^d \arrow \Q^d$  and $A^{tr}_{\psi} : \Q^{d'} \arrow
\Q^{d'}$, hence
 the non-zero spectra of $A^{tr}_\vp $  and 
$A^{tr}_{\psi}$ are the same. \qed

\section{Balanced pairs and Proximality}\label{sec-prox}

The main result of this section is Theorem \ref{homprox}, in which we show
that  if
$\vp$ and $\psi$ are strong Pisot substitutions that satisfy both GCC and the
no cycle condition on periodic words, and $h: \T_\vp \arrow
\T_\psi$ is a homeomorphism, then $h$ carries composants containing proximal
tilings in
$\T_\vp$ to  composants containing proximal tilings in $\T_\psi$.

Suppose that $\vp$ is weak Pisot. Recall that $ \A = \{1, \ldots ,
d\}$ is the alphabet, $A = A_\vp$ the incidence matrix, and $\cL_\vp = \cL
\subseteq \A^*$ the language of $\vp$.  Let
$ \omega_L = (\omega_1, \ldots, \omega_d)$ and
$\omega_R
$ be positive left and right Perron-Frobenius eigenvectors (respectively)
for
$A$.  For $i \in \A$, let $ |i|_g :=
\omega_i$, and for $w = w_1 \ldots w_n \in \A^*$, let $|w|_g := \Sigma_{i
= 1}^n |w_i|_g$, $|w|:= n$ and $l(w) := (a_1, \ldots, a_d)^{tr}$,
where $a_i$ is the number of occurrences of the letter $i$ in $w$ and $v^{tr}$
is the transpose of $v$.  A {\em balanced pair} for
$\vp$ is a pair $\left( \begin{array}{c} u
\\ v
\end{array} \right)$ with $u, v
\in \cL$ and $l(u) = l(v)$.  A {\em geometrically balanced pair} for $\vp$ is
a pair
 $\left( \begin{array}{c} u
\\ v
\end{array} \right)$ with $u, v \in \cL$ and $|u|_g = |v|_g$. Trivially, any
balanced pair is geometrically balanced. If $\vp$ is strong Pisot, then the
entries of
$\omega_L$ are independent over $\Q$, so that if $|u|_g = |v|_g$, then  $l(u)
= l(v)$, and any geometrically balanced pair for $\vp$ is balanced.  As the
Morse-Thue substitution ($1
\arrow 12$,
$2
\arrow 21$) shows, a geometrically balanced pair need not be balanced in
general.

 Define
\[ BP(\vp)
= \left\{
\left( \begin{array}{c} u
\\ v
\end{array} \right) :  \left( \begin{array}{c} u
\\ v
\end{array} \right) \textrm{ is a balanced pair for } \vp\right\} \] and
$$GBP(\vp) = \left\{\left( \begin{array}{c} u
\\ v
\end{array} \right) : \left( \begin{array}{c} u
\\ v
\end{array} \right)\textrm{ is a geometrically balanced pair for }
\vp\right\}.$$     If $x = \left( \begin{array}{c} u
\\ v
\end{array} \right) \in BP(\vp)$, the {\em  dual} of $x$, written
$\overline{x}$, denotes the balanced pair  $\left( \begin{array}{c} v
\\ u
\end{array} \right) $.  If
$\left(
\begin{array}{c} u
\\ v
\end{array} \right) \in BP(\vp)$, and there are $\left( \begin{array}{c} u_1
\\ v_1
\end{array} \right)$,  $\left( \begin{array}{c} u_2
\\ v_2
\end{array} \right) \in BP(\vp)$ so that
$u = u_1 u_2$, $v = v_1 v_2$, then $\left( \begin{array}{c} u
\\ v
\end{array} \right)$ is {\em reducible} and we write $\left( \begin{array}{c}
u
\\ v
\end{array} \right) = \left( \begin{array}{c} u_1
\\ v_1
\end{array} \right) \left( \begin{array}{c} u_2
\\ v_2
\end{array} \right)$.  Otherwise,
  $ \left( \begin{array}{c} u
\\ v
\end{array} \right)$ is {\em irreducible}.  We make similar definitions for
elements of
$GBP(\vp)$.  Note that any balanced pair (geometrically balanced pair,
respectively) factors uniquely as a finite product of  irreducible balanced
pairs (geometrically balanced pairs, respectively).  Let
$\A_{BP}$ be the (possibly infinite) alphabet of irreducible balanced pairs,
and let the substitution $\vp_{BP} : \A_{BP} \arrow (\A_{BP})^*$ be given  by
$\vp_{BP}\left(
\left(
\begin{array}{c} u
\\ v
\end{array} \right) \right) =
\left( \begin{array}{c} \vp(u)
\\ \vp(v)
\end{array} \right)$, factored as a word in $(\A_{BP})^*$.  One can define
$\A_{GBP}$ and
$\vp_{GBP}$ similarly for irreducible geometrically balanced pairs.  An
irreducible balanced pair $ \left( \begin{array}{c} u^0
\\ v^0
\end{array} \right)$  is {\em essential} if for each $n \in \N$, there is $
\left( \begin{array}{c} u^{-n}
\\ v^{-n}
\end{array} \right) \in \A_{BP}$ so that $ \left( \begin{array}{c} u^0
\\ v^0
\end{array} \right)$   is a factor of
$\vp^n_{BP} \left( \left( \begin{array}{c} u^{-n}
\\ v^{-n}
\end{array} \right) \right)$.    An {\em essential geometrically balanced
pair} is defined similarly. Since $\vp$ is primitive,   the trivial pair
$\left( \begin{array}{c} i
\\ i
\end{array} \right)$ is essential for any $i \in A$.

Let $\A_{EBP}$ be the alphabet consisting of essential balanced pairs for
$\vp$. Note that if $\left( \begin{array}{c} {u}
\\ {v}
\end{array} \right)$ is an essential balanced pair, then $\vp_{BP}\left(
\left( \begin{array}{c} {u}
\\ {v}
\end{array} \right) \right)$ is a product of essential factors.  That is, 
$\vp_{BP}$ restricted to
$\A_{EBP}$ determines a substitution
$\vp_{EBP}: \A_{EBP} \arrow (\A_{EBP})^*$.

\begin{example}  \label{esspairs} Essential Balanced Pairs \end{example} 
Define the substitution
$\vp$ as follows:
$$\vp(1) = 11122, \vp(2) = 12.$$  The balanced pairs $\left( \begin{array}{c}
12
\\ 21
\end{array} \right)$, $\left( \begin{array}{c} 112
\\ 211
\end{array} \right)$, and $\left( \begin{array}{c} 1122
\\ 2121
\end{array} \right)$ (and their duals) are all essential:  $$\vp \left(
\left(
\begin{array}{c}12
\\21
\end{array} \right) \right) = \left( \begin{array}{c}1
\\1
\end{array} \right) \left( \begin{array}{c}112
\\211
\end{array} \right) \left( \begin{array}{c}21
\\12
\end{array} \right) \left( \begin{array}{c}2
\\2
\end{array} \right), $$

$$\vp \left( \left( \begin{array}{c}112
\\211
\end{array} \right) \right) = \left( \begin{array}{c}1
\\1
\end{array} \right) \left( \begin{array}{c}112
\\211
\end{array} \right) \left( \begin{array}{c}21
\\12
\end{array} \right) \left( \begin{array}{c}112
\\211
\end{array} \right) \left( \begin{array}{c}21
\\12
\end{array} \right)  \left( \begin{array}{c}2
\\2
\end{array} \right), $$ and 
$$\vp \left( \left( \begin{array}{c}1122
\\2121
\end{array} \right) \right) = \left( \begin{array}{c}1
\\1
\end{array} \right) \left( \begin{array}{c}112
\\211
\end{array} \right) \left( \begin{array}{c}21
\\12
\end{array} \right) \left( \begin{array}{c}1122
\\2121
\end{array} \right) \left( \begin{array}{c}1
\\1
\end{array} \right) \left( \begin{array}{c}21
\\12
\end{array} \right)  \left( \begin{array}{c}2
\\2
\end{array} \right). $$     Let $a, b, c $ (and $\overline{a}, \overline{b},
\overline{c}$) denote the nontrivial essential balanced pairs above in the
order given, and denote the trivial balanced pairs by the
associated letter of $\A_\vp$.  Then $\vp_{EBP}$ must include at least the
information:
$$\vp_{EBP}(a) = 1 b \overline{a} 2, \vp_{EBP}(b) = 1 b \overline{a} b
\overline{a} 2, \vp_{EBP}(c) = 1 b \overline{a} c 1
\overline{a} 2,$$ along with the definitions $\vp_{EBP}(\overline{a}) = 1
\overline{b}
a 2$, etc. 

In Example \ref{rewriting}, we provide, for a
rewriting of
$\vp$, 
an essential geometrically balanced pair that is not balanced. 
\qed

\

Suppose that $T \in \T_\vp$ is fixed by $\Phi$ and  the tile $T_0$ of $T$
containing $0$ is a translate of the prototile $P_a = [0, \omega_a)$, say  $ T_0
= P_a - t_0$.  We will use the location of 0 in $T_0$ to define a {\em
one-cut rewriting of $\vp$ in the letter $a$} according to the two cases
below:

{\em Case 1:}  $t_0 > 0$ (0 is in the interior of $T_0$).  In this case,
$\vp(a) = p a s$ for some (unique) nonempty $p, s$ in $\A^*$ with $ |p|_g <
\lambda t_0$ and $|s|_g < \lambda(w_a - t_0)$.  We modify the alphabet $\A$
by `splitting' the letter $a$ into  two letters $a_1$ and $a_2$.  Let
$\A' = \A \cup \{a_1, a_2\} \setminus \{a\}$, and let $\alpha : \A
\arrow \A'$ be the morphism that takes $a$ to $a_1 a_2$ and $b$ to $b$ if $b
\neq a$.  Define the  substitution
$\vp' :
\A'
\arrow (\A')^*$ as follows: \[ \begin{array}{ll} \vp'(a_1) & = \alpha(p) a_1,
\\
\vp'(a_2) & =  a_2 \alpha(s), \\ \vp'(b) &= \alpha(\vp(b)) \textrm{ if } b
\neq a_1, a_2.
\end{array} \]  Finally, let $\tilde{\vp}$ be the substitution obtained from
$\varphi'$ by rewriting, using starting rule $\{a_2\}$ and stopping rule
$\{a_1\}$.

{\em Case 2:} $t_0 = 0$.  Let the tile type of $T_{-1}$ be $b$.  In this
case, there are  $s, p \in \A^*$ such that $\vp(a) = as$, $\vp(b) = pb$. 
Let $
\tv$ denote the substitution obtained from $\vp$ by rewriting, using starting
rule $\{a\}$ and stopping rule $\{b\}$.

\begin{example} \label{rewriting} One-cut rewriting \end{example} Again,
let 
$\vp$ be given  by:
$\vp(1) = 11122, \vp(2) = 12.$ \\ The  bi-infinite word
$\ldots 11122\ 1\underline{1} 122\ 11122\  12\  12 \ \ldots$ is associated
with a tiling $T$ where 0 occurs in the interior of the tile associated with
the underlined 1 and $T$ is fixed under $\Phi$.

We use the location of $0$ to split $1$ into $1_11_2$.  Then
$\A' =
\{1_1, 1_2, 2\}$, and $\vp'$ is defined by: $$\vp'(1_1) = 1_11_21_1,
\vp'(1_2) = 1_21_11_222,
\vp'(2) = 1_11_22.$$  In rewriting using the starting rule $1_2$ and the
stopping rule $1_1$, one obtains the three words/letters $$ a := 1_21_1,
b := 1_2 2 1_1; c := 1_2 2 2 1_1.$$  The one-cut rewriting $\tv$ is defined
by:
$$\tv(a) := \vp'(1_21_1) = 1_21_11_2221_11_21_1 := aca, \tv(b) =
acba,
\tv(c) = acbba.$$  
The geometrically balanced pair $\left(
\begin{array}{c} ac
\\ bb
\end{array} \right) $ for $\tv$ is not balanced. The substitution $\vp$ is
strong Pisot, with eigenvalues of $ 2 + \sqrt{3}$ and $2 - \sqrt{3} $, while
$\tv$ is weak Pisot, with eigenvalues of $ 2 + \sqrt{3}, 2 - \sqrt{3} $ and
$0$.\qed

\

There is a natural map that takes an arbitrary tiling $T \in
\T_\vp$ to a tiling
$\tT \in \T_{\tv}$, where $\tT$ is obtained by marking $T$ as follows.  In
Case 1,
$T$ is marked at  the cut point in each tile of type $a$.  In Case 2, $T$ is
marked at the beginning of each tile of type $a$ that is preceded by a tile
of type $b$.  Either marking  can be associated naturally with a tiling  in $
\T_{\vp'}$ (Case 1) or
$\T_\vp$ (Case 2).  These tiles are then  amalgamated and relabeled,
according to the rewriting (as seen in Example \ref{rewriting}).  The
correspondence
$T
\mapsto
\tT$ is a homeomorphism that commutes with inflation and substitution
as well as  the flow on both spaces.  That is,
$\widetilde{\Phi(T)} = \tilde{\Phi}(\tT)$  (where
$\td{\Phi}$ denotes inflation and substitution in $\T_{\tv}$) and
$\widetilde{T -t} =
\tT - t$.  

\begin{lemma}If $\tv$ is obtained from $\vp$ by a one-cut rewriting,  then
$\tv$ is proper.  If in addition, $\vp$ is weak Pisot and satisfies the no  cycle
condition, then $\tv$  is weak Pisot.
\end{lemma}

Proof: The first statement follows directly from the definition of one-cut 
rewriting.   Since there is a homeomorphism of $\T_\vp$ and $\T_{\tv}$ that
conjugates $\Phi$ and $\tilde{\Phi}$, the second statement then follows from
Corollary
\ref{cohomcor}.
\qed

\

  Furthermore, if
$\left( \begin{array}{c}
\td{u}
\\ \td{v}
\end{array} \right)$ is a geometrically balanced pair for $\tv$, there is an
associated  pair
$\left( \begin{array}{c} {u}
\\ {v}
\end{array} \right)$ with
$u, v
\in
\cL_\vp$.  Specifically, if $\tv$ is determined by Case 2, then $u $ and $v$
begin with
$a$ and end with
$b$, and $\td{u}$, $\td{v}$ result from rewriting $u, v$.  If $\tv$ is
determined by Case 1, then $\td{u}$ and $\td{v}$ are rewritings of $u' =
a_2x'a_1$ and $v' = a_2y'a_1$ with $x' = \alpha(x)$, $y' = \alpha(y)$ for
some $x, y \in
\cL_\vp$; let $u = ax$, $v = ay$.

\begin{lemma}  If $\tv$ is  a one-cut rewriting of $\vp$, and
$\left( \begin{array}{c}
\td{u}
\\ \td{v}
\end{array} \right) \in GBP(\tv)$, then the associated pair $\left(
\begin{array}{c} {u}
\\ {v}
\end{array} \right) $ is a balanced pair for $\vp$ (not necessarily
irreducible).
\end{lemma}

Proof:  Since $|u|_g = |v|_g$, and $\vp$ is strong Pisot,  $l(u) = l(v)$.
\qed

\

In the following, if $T
\in
\T_\vp$,
$\gamma(T)$ will denote the corresponding strand in
$\T^S_\vp$, and $\Phi$ will denote the inflation and substitution map in
both 
$\T_\vp$ and
$\T^S_\vp$. The flow in $\T_\vp$ will be denoted by $(T, t) \mapsto T - t$,
and we assume that
$w_R$ has been normalized so that $\gamma(T - t) = \gamma(T) - t w_R$.  

  By
a {\em state}
$I$, we mean a line segment of length one, parallel to a coordinate axis and
meeting
$E^s$.  We take
$I$ to be half-open, including its initial vertex but not its terminal
vertex.  If $I$ is a state, $\hat{\Phi}(I)$ will denote the unique state
contained in $\Phi(I)$.  Thus
$\hat{\Phi}$ maps  states to  states. If $\gamma \in \T^S_\vp$ or $\gamma$ is
a finite strand that intersects $E^s$, the state determined by $\gamma $ is
the segment (that is, edge) of
$\gamma$ that meets
$E^s$.

Given tilings $T$ and $T'$ in $\T_\vp$ with strands $\gamma = \gamma(T)$ and
$\gamma' = \gamma(T')$, we say that  $\gamma $ and
$\gamma' $ are {\em coincident at zero} if the states determined by $\gamma $
and
$\gamma' $ are identical;  $\gamma $ and
$\gamma' $ are {\em coincident } if there is $t \in \R$ so that $\gamma  - t
\omega_R$ and
$\gamma'  - t \omega_R$ are coincident at zero.

\

{\parindent=0pt {\bf  Geometric Coincidence Condition (GCC):} If $I$, $J$ }
are states whose vertices are equivalent mod $\Z^d$, then for every $\epsilon
> 0$, there is  $t
\in
\R$ such that
$|t| <
\epsilon$, $I + t \omega_R$ and $J + t \omega_R$ are states, and for some $n
\geq 0$, $\hat{\Phi}^n(I + t \omega_R) = \hat{\Phi}^n(J + t \omega_R)$
(equivalently, for some
$n \geq 0$,  ${\Phi}^n(I + t \omega_R) $ and ${\Phi}^n(J + t \omega_R)$ are
coincident at 0).

\

 If $\vp$ is strong Pisot, and $\gamma \in \T^S_\vp $, let $v(\gamma)$ be any
vertex of $\gamma$. Define
$g: \T^S_\vp \arrow \cT^d$  by $ g(\gamma) = v(\gamma)
\ (\textrm{mod }
\Z^d$).  Then $g$ is a continuous surjection and $g(\gamma - t \omega_R) =
g(\gamma) - t \omega_R \ (\textrm{mod }
\Z^d)$ for all $t \in \R$.  Also, if $F_A: \cT^d \arrow \cT^d$ is defined by
$F_A(p) = A p \  (\textrm{mod }
\Z^d)$, then $g \circ \Phi = F_A \circ g$.  In the case $\vp$ is unimodular,
the map $g$ is called {\em geometric realization}.  

If
$\vp$ is not unimodular, define $\hat{g} : \T^S_\vp \arrow \inv F_A$ by
$$\hat{g}(\gamma) = (g(\gamma), g(\Phi^{-1}(\gamma)), g(\Phi^{-2}(\gamma)),
\ldots );$$ in this case,
$\hat{g}$ is geometric realization.  Let $(\underline{z})_t$ denote the flow
on $\inv F_A$ defined by $$(z_0, z_1, \ldots )_t = (z_0 - t\omega_R, z_1
- (t/\lambda)
\omega_R, z_2 - (t/\lambda^2)
\omega_R, \ldots),$$ where all coordinates are taken mod $\Z^d$,  and let
$\hat{F}_A$ be the shift homeomorphism on
$\inv F_A$ given by $$\hat{F}_A (z_0, z_1, \ldots ) = (F_A(z_0), z_0, z_1,
\ldots).$$ Then
$\hat{g}
\circ
\Phi =
\hat{F_A}
\circ
\hat{g}$, and
$\hat{g}(\gamma - t
\omega_R) = (\hat{g}(\gamma))_t
$ for all $t \in \R$.

Note that if $\vp$ is unimodular, then $F_A: \cT^d \arrow \cT^d$ is a
homeomorphism, so that projection onto the first coordinate, $\pi_0$, yields
a homeomorphism of $\inv F_A$ with $\cT^d$ that conjugates the $\R$- and
$\Z$-actions on $\inv F_A$ with those on $\cT^d$ and produces the commuting
diagram:

{\parindent=0pt
\begin{picture}(25, 80)(-130, -40)

\put(80, -8){$\pi_0$}
\put(34, -32){$\cT^d$}
\put(-10, -8){$g$}

\put(-20, 20){$ \T_\vp^S$}
\put(70, 20){$\inv F_A $}
\put(-10,10){\vector(1,-1){30}}
\put(-0,25){\vector(1,0){65}}
\put(35, 30){$\hat{g}$}

\put(85,10){\vector(-1,-1){30}}

\end{picture}  

For} efficiency in the following, we will denote geometric realization in the
unimodular case by $\hat{g}$ as well.

For a strong Pisot substitution $\vp$, the {\em coincidence
rank of
$\vp$} is defined as 
\begin{eqnarray*}
cr_\vp = &\max \{n: \textrm{ there are } \gamma_1,
\ldots, \gamma_n \in
\T_\vp^S \textrm{ and } \underline{x} \in \inv F_A  \textrm{ such that } \\ &
 \gamma_i
\in
\hat{g}^{-1}(\underline{x}) \textrm{ for }  1 \leq i \leq n \textrm{ and } 
 \gamma_i, \gamma_j  \textrm{ are not coincident if } i \neq j\}.
\end{eqnarray*}
 A balanced pair $B$ for $\vp$  {\em terminates
with coincidence} provided there is a finite collection $\{B_1, \ldots,
B_k\}$ of irreducible balanced pairs so that (1) for each $n \in \N \cup
\{0\}$, the pair
$\vp^n_{BP}(B)$ factors as a product of the elements of $\{B_1, \ldots,
B_k\}$, and (2) for each $i \in \{1, \ldots, k\}$, there is  $n \in \N$ so
that 
$\vp^n_{BP}(B_i)$ has a trivial balanced pair factor.

Proofs of the following results in the  case $\vp$ is unimodular can be found
in
\cite{mjarek}.  Straightforward modifications of those proofs yield the
non-unimodular case.

\begin{theorem} \label{GCCeq} Suppose that $\vp$ is strong Pisot. Then
$\hat{g}$ is bounded-to-one and $\#\hat{g}^{-1}(\underline{x}) = cr_\vp$
for (Haar) almost every $\underline{x} \in \inv F_A$.  Moreover, the
following are equivalent:
\\(i) $\vp$ satisfies GCC, 
\\(ii) $ cr_\vp = 1$,
\\(iii) $\# \hat{g}^{-1}(\underline{x}) =1$ for (Haar)
almost every $\underline{x} \in \inv F_A$,
\\(iv)  every balanced pair $B$ for $\vp$ terminates with coincidence,
\\(v)  there are $a, b \in \A$ with $a \neq b $ so that the balanced pair
$\larr a
\ b
\\ b
\ a
\rarr$ terminates with coincidence,
\\(vi) the tiling flow, $(T, t) \mapsto T - t$ has pure discrete spectrum.
\\If any of (i)-(vi) holds, then $\{\gamma \in \T_\vp^S : 
\#\hat{g}^{-1}(g(\gamma)) = 1\}$ is a set of full measure (with respect to
the unique ergodic flow-invariant measure on $\T_\vp^S$) that contains a
dense
$G_\delta$.
\end{theorem}

Tilings $T, T' \in \T_\vp$ are said to be {\em forward proximal} if there is
a sequence $\{t_k\}$ of real numbers  so that
$t_k \arrow \infty$ and $\lim_{k \arrow \infty} d(T - t_k, T'-t_k) = 0$. 
If there is a sequence  $\{t_k\}$ of real numbers so that
$t_k \arrow -\infty$ and $\lim_{k \arrow \infty} d(T - t_k, T'-t_k) = 0$,
then
$T$ and $ T'$ are said to be {\em backward proximal}.  If $T$ and $T'$ are
either forward or backward proximal, they are {\em proximal}. 
Note
that if $T$ and $T'$ are proximal, so are: $T- t$ and $T' - t$ for each $t
\in
\R$; $\Phi(T)$ and $\Phi(T')$; and $\Phi^{-1}(T)$ and $\Phi^{-1}(T')$ (since
$\Phi^{-1}(T - t) = \Phi^{-1}(T) - t/\lambda$).

\begin{example} \label{proximal} Proximality \end{example}  Let 
$\vp$ be defined by
$\vp(1) = 11122, \vp(2) = 12.$ \\ The following three bi-infinite words
 (given in pairs) represent tilings that are fixed under $\Phi$.  In one
case, the origin of the tiling is at an endpoint of a tile and the
corresponding word is fixed under $\vp$.  The origin of each of the
remaining two tilings  is in the interior of the tile, and the associated
words can be obtained by iterating around a fixed point 
associated with the underlined symbol.  All three are (pairwise)  proximal
in both directions.  Spacing  is used to indicate the balanced pair
structure for each pair of words.
$$ \begin{array}{c} \ldots  1\ 112\ 2 1\ 1\underline{1}2\ 2 1\ 2\ 1\  211
\ldots\\ \ldots  1\ 2 11\ 12\ 2 1\underline{1}\ 12\ 2\ 1\ 112
\ldots \end{array} $$
$$ \begin{array}{c} \ldots  1\ 112\ 2 1\ 1\underline{1}2 2\ 1\  2 1\  2\ 1\ 1
\ldots\\ \ldots  1\ 2 11\ 12\ 2 12  .  1\ 1\ 1 2\  2\ 1\ 1 
\ldots \end{array} $$
$$ \begin{array}{c}  \ldots  1\ 2 11\ 12\ 2\ 1\ \underline{1} 12\ 2 1\ 112
\ldots\\ \ldots  1\ 2 11\ 12\ 2\ 1\ 2  .  1 1\ 1 2\ 211 
\ldots \end{array} $$  Note that each pair could also be generated by
iterating  under
$\vp_{BP}$ (or $\vp_{EBP}$) the balanced pair in which the fixed points of
the bi-infinite words appear. (See Example \ref{esspairs} for $\vp_{EBP}$.) 
\qed

\begin{prop}\label{proxgeom} Suppose that $\vp$ is strong Pisot, and let $T, T'
\in \T_\vp$.  If $T$ and $T'$ are
either forward or backward proximal, then $\gamma(T)$ and $\gamma(T')$ have
the same geometric realization. If $\vp$ satisfies GCC and
$\gamma(T) $ and $\gamma(T')$ have the same geometric realization, then $T$
and $T'$ are proximal in both directions.  
\end{prop}

Proof:  Suppose that $T$ and $T'$ are either forward or backward proximal.
Then so are the images of $\gamma(T)$ and $\gamma(T')$ under the irrational
flow on $\cT^d$ (or $\inv F_A$, in the non-unimodular case). Since 
irrational flow is an isometry,  this implies that
$g(\gamma(T)) = g(\gamma(T'))$  (or $\hat{g}(\gamma(T)) =
\hat{g}(\gamma(T'))$, respectively).

Suppose that $\vp$ satisfies GCC, and that $\gamma= \gamma(T)$ and
$\gamma' = \gamma(T')$ have the same geometric realization.   Choose $n_k
\arrow \infty$ so that 
$\Phi^{-n_k}(\gamma)
$ and
$\Phi^{-n_k}(\gamma') $ converge, say to $\eta$ and $\eta'$ respectively. 
Then $\eta$ and $\eta'$ have the same geometric realization and, using GCC,
there are $t > 0$ and $n_0 \in \N$ so that if $I$ and $I'$ are the states
determined by $\eta + t \omega_R$ and $\eta'  + t \omega_R$, then
$\hat{\Phi}^{n_0}(I) =
\hat{\Phi}^{n_0}(I')$.  Moreover, we may choose $t$ so that neither $I$ nor $I'$
has a vertex on $E^s$.  Then, for sufficiently large $k$, the states $I_{n_k}$
and $I'_{n_k}$ determined by ${\Phi}^{-n_k}(\gamma) + t \omega_R$  and
${\Phi}^{-n_k}(\gamma') + t \omega_R$ satisfy:  $I = I_{n_k} + u_k$ and $I' =
I'_{n_k} + u_k$ with $u_k \in \R^d$ and $u_k \arrow 0$ as $k \arrow
\infty$.  Thus $\hat{\Phi}^{n_0}(I_{n_k} + u_k)$  =
$\hat{\Phi}^{n_0}(I'_{n_k} + u_k)$ for sufficiently large $k$.  That is,
$\Phi^{n_0}  ({\Phi}^{-n_k}(\gamma))$ and $\Phi^{n_0} 
({\Phi}^{-n_k}(\gamma'))$ share an edge, call it $L_k$, whose initial vertex,
$\min L_k$, satisfies $\langle \min L_k, \omega_R\rangle < 0$, provided $k$ is
large enough so that  $\langle t\omega_R + u_k, \omega_R\rangle 
> 0$ (here $\langle, \rangle$ is the usual Euclidean inner product).  It follows
that $\gamma$ and $\gamma'$ coincide along ${\Phi}^{n_k - n_0}(L_k) $.  Since
$\langle A^{n_k - n_0} (\min L_k), \omega_R \rangle$ goes to $- \infty$ as $k
\arrow \infty$ and the length of ${\Phi}^{n_k - n_0}(L_k) $ goes to $\infty$ as
$k \arrow \infty$, $T$ and $T'$ are backward proximal.  Choosing $t$ as above but
with $t < 0$ shows that $T$ and $T'$ are forward proximal.
 \qed

\begin{corollary} \label{proxequiv} Suppose that  $\vp$ is a strong Pisot
substitution that satisfies GCC. Then proximality is an equivalence
relation, and if $T$ and
$T'$ are either forward or backward proximal, they are proximal in both
directions.  The proximality equivalence class of $T$ is exactly the
collection of
$T'$ for which
$\gamma(T)$ and $\gamma(T')$ have the same geometric realization.
\end{corollary}

  Suppose that
$T$ and $T'$ in $\T_\vp$  are proximal, and that $\vp$ is a strong Pisot
substitution satisying GCC.  If $\gamma = \gamma(T)$ and
$\gamma' = \gamma(T')$, there are times $t_k$ (for $k \in \Z$) such
that $\lim_{k
\arrow -\infty}t_k = -\infty$,
 $\lim_{k
\arrow \infty}t_k = \infty$, and    $\gamma - t_k \omega_R$ and $\gamma' -
t_k
\omega_R$ are coincident at zero for all $k$.  Let  $$NC = \{t: \gamma - t
\omega_R
\textrm{ and } \gamma' - t
\omega_R \textrm{ are {\em not} coincident at zero} \}.$$  If $D$ is a
component of
$NC$, then $D$ is a bounded interval of the form $D = [a, b)$ (recall that we
take the segments defining states  to be closed on the left and open on the
right).  Let
$\gamma_{[a, b)} = \gamma \cap ( \cup_{t \in [a, b)} (E^s + t\omega_R))$, and
let $u$ be the word in $\cL_\vp$ determined by $\gamma_{[a, b)}$.  Similarly,
let $v$ be the word determined by   $\gamma'_{[a, b)}$.  Then  $\left(
\begin{array}{c} {u}
\\ {v}
\end{array} \right)$ is a geometrically balanced, hence balanced pair for
$\vp$.  If  $\gamma_{[a, b)}$ and
$\gamma'_{[a, b)}$ do not intersect in their interiors,  $\left(
\begin{array}{c} {u}
\\ {v}
\end{array} \right)$ is irreducible.  In this case we say that
$\left( \begin{array}{c} {u}
\\ {v}
\end{array} \right)$ is {\em obtained from a bubble in a proximal pair}.

\begin{prop} \label{bubble} Suppose that $\vp$ is strong Pisot and satisfies
GCC. Let $\left( \begin{array}{c} {u}
\\ {v}
\end{array} \right)$ be an irreducible balanced pair for $\vp$. Then
$\left( \begin{array}{c} {u}
\\ {v}
\end{array} \right)$ is an essential balanced pair for $\vp$ if and only if 
$\left( \begin{array}{c} {u}
\\ {v}
\end{array} \right)$ is obtained from a bubble in a proximal pair or is a
trivial pair $\left( \begin{array}{c} {i}
\\ {i}
\end{array} \right)$.
\end{prop}

Proof:  Suppose that  $\left( \begin{array}{c} {u}
\\ {v}
\end{array} \right)$ is obtained from a bubble in the proximal pair $T, T'$. 
For each $n \in \N$,
$\Phi^{-n}(T)$ and $\Phi^{-n}(T')$ are also proximal.  There is then a
`geometrical bubble' in $\gamma(\Phi^{-n}(T))$  and $\gamma(\Phi^{-n}(T'))$
that maps, under
$\Phi^n$, over the geometrical bubble that determines $\left(
\begin{array}{c} {u}
\\ {v}
\end{array} \right)$. Thus there is an irreducible balanced pair
$\left( \begin{array}{c} {u^{-n}}
\\ {v^{-n}}
\end{array} \right)$ (the linguistic equivalent of the geometrical bubble in
$\gamma(\Phi^{-n}(T))$  and $\gamma(\Phi^{-n}(T')))$ that maps, under
$\vp^n_{BP}$, into a word having $\left( \begin{array}{c} {u}
\\ {v}
\end{array} \right)$ as a factor, thus $\left( \begin{array}{c} {u}
\\ {v}
\end{array} \right)$ is essential.

Conversely, if $\left( \begin{array}{c} {u}
\\ {v}
\end{array} \right)$ is essential, for each $n \in \N$, let $\left(
\begin{array}{c} {u^{-n}}
\\ {v^{-n}}
\end{array} \right)$ be an irreducible balanced pair with $\left(
\begin{array}{c} {u}
\\ {v}
\end{array} \right)$  a factor of $\vp^n_{BP}\left( \left( \begin{array}{c}
{u^{-n}}
\\ {v^{-n}}
\end{array} \right) \right)$.  For each $n$, let $\gamma_{-n}, \gamma'_{-n}$
be a pair of finite strands that  realize the patterns $u^{-n}$ and $v^{-n}$
and have the same initial and terminal points (so that $\gamma_{-n}
\cup\gamma'_{-n}$ is a geometrical bubble) located a bounded (minimum)
distance from the origin in $\R^d$.  By adjusting by translation in the $E^u$
direction, we can insure that  $\Phi^n(\gamma_{-n}) \cup
\Phi^n(\gamma'_{-n})$ contains a bubble linguistically equivalent to $\left(
\begin{array}{c} {u}
\\ {v}
\end{array} \right)$ that meets $E^s$.  There is a subsequence $\{n_i\}$ so
that the increasingly long strands $\Phi^{n_i}(\gamma_{-n_i}) $ converge to a
strand $\gamma \in \T^S_\vp
$.  There is a further subsequence so that the strands
$\Phi^{n_{i_j}}(\gamma'_{-n_{i_j}})
$ converge to a strand $\gamma' \in \T^S_\vp
$.    It is clear that $\gamma$ and
$\gamma'$  have the same geometric realization, so that if  $T$ and $T'$
satisfy
 $\gamma = \gamma(T)$ and $ \gamma' = \gamma'(T')$, then $T$ and $T'$ are
proximal.   Finally, $\left( \begin{array}{c} {u}
\\ {v}
\end{array} \right)$ is obtained from a bubble in the proximal pair $T, T'$.
\qed

\

Recall that $\A_{EBP}$ is the alphabet consisting of essential balanced pairs
for
$\vp$.  According to the following result, if   $\vp$ is strong
Pisot and satisfies GCC, the substitution
$\vp_{EBP}: \A_{EBP} \arrow (\A_{EBP})^*$ is on a finite alphabet. 
\begin{lemma} \label{essfinite} If $\vp$ is strong Pisot and satisfies GCC,
then
$\A_{EBP}$ is finite.
\end{lemma} Proof:  If $\A_{EBP}$ is infinite, then according to Proposition
\ref{bubble}, there are
$T^n, S^n
\in
\T_\vp$ so that
$T^n$ and
$S^n$ are proximal and $\gamma(T^n)$, $\gamma(S^n)$ determine an irreducible
geometrical bubble of length at least $2n|\omega_R|$. By translating, we may
assume that these geometric bubbles extend at least from $E^s - n\omega_R$
to $E^s + n\omega_R$. Choose a subsequence $\{n_i\}$ so that $T^{n_i}
\arrow T \in
\T_\vp$ and
$S^{n_i} \arrow S \in \T_\vp$.  Since $T^{n_i}$ has the same geometric
realization as $S^{n_i}$ for each $i$, $T$ and $S$ have the same geometric
realization, hence $T$ and $S$ must be proximal.  But $\gamma(T)$ and
$\gamma(S)$ do not intersect, so $T$ and $S$ cannot be proximal.  \qed

\

 If we order
the elements of  $\A_{EBP}$ as $\{e_1, \ldots, e_n\}$ where for $i = 1,
\ldots, d$,
$e_i
$ denotes the trivial balanced pair $\left( \begin{array}{c} {i}
\\ {i}
\end{array} \right)$, then the matrix for $\vp_{EBP}$,
$A_{\vp_{EBP}}$, has the block triangular form $\left( \begin{array}{cc}
A_\vp &  B
\\  0  & C
\end{array}
\right)$. (In particular, $\vp_{EBP}$ is not primitive.)  

For the tiling space $\T_{\vp_{EBP}}$, let the length  of the prototile
corresponding to balanced pair $e_i = \larr u_i \\ v_i \rarr$ be given by
$|u_i|_g = \omega_i$, let $|e_{i_1} \ldots e_{i_k}|_g = \Sigma_{s = 1}^t
w_{i_s}$ be the geometrical length of a word in $A^*_{\vp_{EBP}}$, and let
$\Phi_{EBP}$ denote the inflation and substitution homeomorphism on
$\T_{EBP}$.  (Ordinarily, the invertibility  of a substitution $\vp$ or
associated map on strand space $\Phi$ is {\em recognizability}, which depends
on primitivity.  But in this case, $\Phi$ is invertible, so 
$\Phi_{EBP}$ is invertible.)  Then $|\Phi_{EBP}(e_i)|_g = \lambda |e_i|_g$
for $i = 1, \ldots, m$, so that $\omega_{L, \vp_{EBP}} := (\omega_1,
\ldots,
\omega_m)$ is a positive left eigenvector for $A_{\varphi_{EBP}}$ with
eigenvalue $\lambda$.

From the primitivity of $\vp$, it follows that the positive
right Perron-Frobenius eigenvector $\omega_R =  (f_1
\ldots, f_d)^{tr} $ of $A_\vp$, normalized so that $|\omega_R|_1 = \Sigma_{i
= 1}^d |f_i| = 1$, has entries 
 $f_i$ equal to the  frequency of occurrence of the tiles of type $i$ in
any tiling $T =
\{T_n\}^\infty_{n = -\infty} \in \T_\vp$: for $i = 1, \ldots, d$,
$$f_i = \lim_{N \arrow
\infty}
\frac{1}{2N+1}
\#\{ n: -N \leq n \leq N \textrm{ and } T_n \textrm{ has type } i\} 
.$$  We extend this result to $\T_{\vp_{EBP}}$.

\begin{lemma} \label{freq} Assume that $\vp$ is strong Pisot and satisfies
GCC, let
$f_i$ be as above for $i = 1, \ldots, d$, and for $i = d+ 1, \ldots, m =
|\A_{\vp_{EBP}}|$, let $f_i = 0$.  Then, for any tiling $T =
\{T_n\}^\infty_{n = -\infty} \in \T_{\vp_{EBP}}$, $$f_i = \lim_{N \arrow
\infty}
\frac{1}{2N+1}
\#\{ n: -N \leq n \leq N \textrm{ and } T_n \textrm{ has type } e_i\} 
$$ for all  $i = 1, \ldots, m$, and the limit is uniform in $T$.
  \end{lemma}

Proof:  The fact that
$\vp$ is primitive and satisfies  GCC implies that  $(A_{\vp_{EBP}})^n $ has
the form $
\left(
\begin{array}{cc} (A_\vp)^n &  B_n \\  0  & C^n \end{array}
\right)$ with $B_n$  strictly positive for large enough $n$, say $ n \geq
M$. Then $$A^{lM}_{\vp_{EBP}}  = \left( \begin{array}{cc}
      A^{lM}_\vp &  [\Sigma_{i = 1}^l
A_\vp^{(l-i)M} B_M C^{(i-1)M}] \\ 0 & C^{lM}
\end{array} \right) $$ for $l = 1, 2, \ldots$. Let 
$\omega_L = (\omega_1, \ldots,
\omega_d)$ and $v = (\omega_{d+1}, \ldots,
\omega_m)$ so that $(\omega_L, v) A_{\varphi_{EBP}} = (\lambda \omega_L,
\lambda v)$.  Then $$\omega_L (\Sigma_{i = 1}^l
A_\vp^{(l-i)M} B_M C^{(i-1)M}) + vC^{lM} = \lambda^{lm}v.$$
That is, 
 $$ (\Sigma_{i = 1}^l \lambda^{(l-i)M}
\omega_L B_M C^{(i-1)M}) + vC^{lM} = \lambda^{lm}v, \ \ (*)$$ with
$\omega_LB_M$ and
$v$ strictly positive.

Let $\beta$ be the spectral radius of $C$.  Since $C \geq 0$, $\beta$ is a
real eigenvalue of $C$ and since $\omega_LB_M, v$ are strictly positive,
there is a constant $K > 0$ so that $|\omega_L B_M C^{(i-1)M}|_1 \geq
K\beta^{(i - 1)M}$ and $|v C^{lM})|_1 \geq
K\beta^{lM}$ for all $l \in \N$.  Thus, from ($*$), $$\lambda^{lm}|v|_1 \geq
(\beta^{lm} + \Sigma_{i = 1}^l \lambda^{(l-i)M}
\beta^{(i-1)M})K$$ for all $l \in \N$.   This implies that $\beta <
\lambda$.  It follows  (again using the fact that $B_n > 0$ for $n \geq
M$) that, up to scale,  $\omega_{R, \vp_{EBP}} := (\omega_R, 0) =  (f_1,
\ldots, f_d, 0, \ldots 0)$ is the unique eigenvector for $A_{\vp_{EBP}}$
with all entries non-negative and that 
$$ 
\frac{A^{n}_{\vp_{EBP}}w}{|A^{n}_{\vp_{EBP}} w | }
\arrow  \omega_{R, \vp_{EBP}}  
$$  for any non-negative $w \in \R^m$ with at least one nonzero entry.  In
particular, for any $e_i \in \A_{\vp_{EBP}}$ and $k \in \{1, \ldots,
m\}$, $$
\frac{(l(\vp^n_{EBP}(e_i)))_k}{|l(\vp^n_{EBP}(e_i))| }  \arrow f_k , $$
where the numerator represents the $k^{th}$ component of the abelianization
vector $l$ of the word $\vp^n_{EBP}(e_i)$.  Thus, for each $k \in \{1,
\ldots, m\}$ and $ j \in \N$, there is $\epsilon_j > 0$ so that for all $i
\in \{1,
\ldots, m\}$, $$
|\frac{(l(\vp^j_{EBP}(e_i)))_k}{|l(\vp^j_{EBP}(e_i))| }  - f_k | < \epsilon_j
$$ and  $\epsilon_j \arrow 0$ as $j \arrow \infty$.

Now for $N \in \N$, let $w = w_{-N} \ldots w_N$ be the word in
$\A^*_{\vp_{EBP}}$ corresponding to the central portion of a tiling $T =
\{T_n\}_{n = - \infty}^{\infty}$ (that is, $T_n$ has type $w_n$).  For any $j
\in \N$, $w$ can be factored as $w = u \vp^j_{EBP}(e_{i_1}) \ldots
\vp^j_{EBP}(e_{i_t})v$ with $|u|, |v|$ bounded independently of $T$ and $N$. 
It follows that 
$$  \frac{ -\epsilon_j - \frac{|u| + |v|  }{ \Sigma_{s =
1}^t|l(\vp^j_{EBP}(e_{i_s}))|} f_k  }
{ \frac{|u| + |v|  }{ \Sigma_{s =
1}^t|l(\vp^j_{EBP}(e_{i_s}))| } + 1  }  
$$ 
$$ \leq  \frac{
\frac{\Sigma_{s=1}^t(l(\vp^j_{EBP}(e_{i_s})))_k}{
\Sigma_{s=1}^t|l(\vp^j_{EBP}(e_{i_s}))| }  - f_k -
\frac{|u| + |v|  }{ \Sigma_{s = 1}^t|l(\vp^j_{EBP}(e_{i_s}))|} f_k  } {
\frac{|u| + |v|  }{ \Sigma_{s = 1}^t|l(\vp^j_{EBP}(e_{i_s}))| } + 1  }   $$
$$ \leq  \frac{\Sigma_{s=1}^t(l(\vp^j_{EBP}(e_{i_s})))_k}
{|u| + |v|  + \Sigma_{s = 1}^t|l(\vp^j_{EBP}(e_{i_s}))|  } -f_k $$
$$ \leq 
\frac{1}{2N+1}
\#\{ n: -N \leq n \leq N \textrm{ and } T_n \textrm{ has type } e_k\} -
f_k$$ 
$$ \leq  \frac{|u| + |v|  + \Sigma_{s=1}^t(l(\vp^j_{EBP}(e_{i_s})))_k}
{ \Sigma_{s = 1}^t|l(\vp^j_{EBP}(e_{i_s}))|  } -f_k $$
$$ \leq  \frac{|u| + |v|  }
{ \Sigma_{s = 1}^t|l(\vp^j_{EBP}(e_{i_s}))|  } + \epsilon_j.$$  Since
$\Sigma_{s = 1}^t|l(\vp^j_{EBP}(e_{i_s}))| \arrow \infty$ as $N \arrow
\infty$, we have $$\limsup_{N \arrow \infty} | \frac{1}{2N+1}
\#\{ n: -N \leq n \leq N \textrm{ and } T_n \textrm{ has type } e_k\} -
f_k | \leq \epsilon_j.$$  Since $\epsilon_j \arrow 0 $ as $j \arrow \infty$,
we have the desired result.
\qed

\

 We have just proved that in the tiling space $\T_{\vp_{EBP}}$, the tile
types corresponding to nontrivial essential balanced pairs occur with zero
frequency in any tiling.  Equivalently, if
$T, T'$ is a proximal pair in $\T_\vp$, then $$\mu \{t : t_0 \leq t \leq t_0
+
\tau, t
\in NC\}/ \tau \arrow 0,$$  uniformly in $t_0$, as $\tau \arrow \infty$
(where $\mu$ is Lebesgue measure, and $NC$ is the previously defined set of
non-coincident times).

\begin{lemma} \label{essgfinite} Suppose that $\vp$ is strong Pisot and
satisfies both GCC and the no cycle condition.  If
$\tv$ is derived from $\vp$ by a one-cut rewriting, then $$\A_{\tv, EGBP}
:=
\left\{\left(
\begin{array}{c} \td{u}
\\ \td{v}
\end{array} \right) : \left( \begin{array}{c} \td{u}
\\ \td{v}
\end{array} \right)\textrm{ is an essential  geom. bal. pair for } \tv
\right\}$$  is finite.
\end{lemma}

Proof:  Suppose that $\left(
\begin{array}{c} \td{u}^{-n}
\\ \td{v}^{-n}
\end{array} \right)$ is a sequence of irreducible geometrically balanced pairs for
$\tv$ with $\left(
\begin{array}{c} \td{u}
\\ \td{v}
\end{array} \right)$ a factor of $\tv^n_{GBP} \left( \left(
\begin{array}{c} \td{u}^{-n}
\\ \td{v}^{-n}
\end{array} \right) \right)$ for $n  \in \N$.  Then the corresponding
balanced pairs $\left( \begin{array}{c} {u}
\\ {v}
\end{array} \right)$ and $ \left( \begin{array}{c} {u}^{-n}
\\ {v}^{-n}
\end{array} \right)$ for $\vp$ have the property that $\left(
\begin{array}{c} {u}
\\ {v}
\end{array} \right)$ is a factor of  $\vp^n_{BP} \left( \left(
\begin{array}{c} {u}^{-n}
\\ {v}^{-n}
\end{array} \right) \right)$ for $n \in \N$.  It follows that the irreducible
factors of $\left( \begin{array}{c} {u}
\\ {v}
\end{array} \right)$ are all essential.  Now, the cut letter $a$ (or the word
$ba$, in the case $\tv$ is obtained from $\vp$ as in Case 2 of the definition
of $\tv$) occurs with bounded gap.  That is, there is an $N \in \N$ so that
if $W \in \cL_\vp$ has length at least $N$, then $W$ contains $a$ ($ba$,
respectively) as a factor.  It follows from Lemma \ref{freq} that there is a
$K
\in
\N$ so that if
$W \in \cL_{\vp_{EBP}}$ with $|W| \geq K$, then $W$ contains a factor that is
itself a product of $N$ trivial balanced pairs.  Thus, if $\left(
\begin{array}{c} \td{u}
\\ \td{v}
\end{array} \right)$ is an essential geometrically balanced pair for $\tv$
that is long enough so that the corresponding  $\left( \begin{array}{c} {u}
\\ {v}
\end{array} \right)$ factors into a product of at least $K$ essential
balanced pairs for $\vp$, then $\left( \begin{array}{c} {u}
\\ {v}
\end{array} \right) = \ldots \left( \begin{array}{c} {i_1}
\\ {i_1}
\end{array} \right)\left( \begin{array}{c} {i_2}
\\ {i_2}
\end{array} \right) \ldots \left( \begin{array}{c} {i_N}
\\ {i_N}
\end{array} \right) \ldots$ with $i_j \in \A$ for $j = 1, \ldots, N$.  There
must be
$j$ with
$i_j = a$ (or
$i_j i_{j+1} = ba$, respectively), hence $\left(
\begin{array}{c} \td{u}
\\ \td{v}
\end{array} \right)$ is not irreducible and thus not essential. It follows
that there are only finitely many essential  geometrically balanced pairs for
$\tv$. 
\qed

\

Suppose, for the remainder of this section, that {\bf  $\vp$ is strong Pisot
and satisfies GCC and the no  cycle condition on periodic words}.  In
addition, {\bf
$\tv$ is obtained from $\vp$ by a one-cut rewriting}.  Let $$\T^P_\vp =
\left\{ \left( \begin{array}{c} T \\ T' \end{array}
\right)  : T, T' \in \T_\vp; T, T' \textrm{are proximal}  \right\} $$ have
the natural (product) topology.  It follows from Proposition \ref{bubble} and
Lemma \ref{essfinite} that $\T^P_\vp  \simeq \T_{\vp_{EBP}}$.  Let
$\tv_{EGBP} :
\A_{\tv, EGBP} \arrow (\A_{\tv, EGBP} )^*$ be the substitution (on a finite
alphabet, by Lemma
\ref{essgfinite}) given by $$\tv_{EGBP}\left( \begin{array}{c}
\td{u} \\ \td{v} \end{array}
\right) = \left( \begin{array}{c}
\td{\vp}  (\td{u}) \\ \td{\vp}(\td{v})  \end{array}
\right),$$ factored as a product of essential geometrical balanced pairs. 
Let
$\A_{\tv} = \{ 1, \ldots,
\td{d}\}$. (Note that the symbol $i \in \A_\vp$ is not equal to the symbol
$i \in \A_{\tv}$, since $\tv$ is a one-cut rewriting of $\vp$ and not an
extension of $\vp$.) If we order the elements of
$\A_{\tv, EGBP} $ as $\{\td{e_1}, \ldots, \td{e}_{n} \}$ where for  $i =
1,
\ldots, \td{d}$,   $\td{e}_i$ denotes the trivial geometrical balanced
pair $
\left(
\begin{array}{c}
i \\ i \end{array}
\right)$ for $\tv$, then the matrix for $\tv_{EGBP}$ has the form
$\left( \begin{array}{cc} \td{A} &  \td{B} \\  0  & \td{C} \end{array}
\right)$ where $ \td{A} = A_{\tv}
$.  Again,  if $\left( \begin{array}{cc} \td{A} &  \td{B} \\  0  & \td{C}
\end{array}
\right)^n  = \left( \begin{array}{cc} \td{A}^n &  \td{B}_n \\  0  & \td{C}^n
\end{array}
\right)$, then $\td{B}_n$ is strictly positive for sufficiently large $n$. 
Thus (see  Lemma \ref{freq}) if $\td{S}  =
\{\td{S_k}\}^\infty_{k = - \infty}
\in
\T_{\tv_{EGBP}}$, the tiles  $\td{S_k}$ are predominantly of types
$\{\td{e_1},
\ldots, \td{e}_{\tilde{d}} \}$.  In fact, there must be $\{k_m\}^\infty_{m =
-
\infty}$ with $\lim_{m \arrow - \infty} k_m = -\infty$, $\lim_{m \arrow 
\infty} k_m = \infty$, and   $\td{S}_{k_{m}}$,  $\td{S}_{k_m + 1}$, \ldots,
$\td{S}_{k_m + |m|}$  all of trivial type for each $m$.  We can interpret
$\td{S}$ as a pair $ \left( \begin{array}{c}
\td{T} \\ \td{T}' \end{array}
\right)$ of tilings in $\T_{\tv}$: if $\td{S_k}$ is a tile of type $ \td{e}_i
=
\left(
\begin{array}{c}
\td{u} \\ \td{v} \end{array}
\right)$, tiling the interval $[a, b)$, then in $\td{T}$ and $\td{T}'$, the
interval
$[a, b)$ is tiled following the patterns $\td{u}$ and $\td{v}$,
respectively.  Under this interpretation,
$\td{T}$ and
$\td{T}'$ are proximal.  Conversely, if $\td{T}$ and $\td{T}'$ are proximal
in
$\T_{\tv}$, then the corresponding $T, T'$ in $\T_\vp$ are proximal and the
pair $
\left(
\begin{array}{c} T \\ T' \end{array}
\right)$ determines $S \in \T_{\vp_{EBP}}$.  The tiling $S$ in turn uniquely
determines
$\td{S} \in \T_{\tv_{EGBP}}$  (for instance, any two consecutive
occurrences of  $\left( \begin{array}{c} a
\\ a
\end{array} \right)$ in $S$ (or  $\left( \begin{array}{c} b
\\ b
\end{array} \right)$ $\left( \begin{array}{c} a
\\ a
\end{array} \right)$ in the case $\tv$ is constructed as in Case 2) uniquely
determines the decomposition of the associated section of  $\left(
\begin{array}{c} \tilde{T} 
\\ \tilde{T'}
\end{array} \right)$ into essential geometrically balanced pairs).

Thus $ \T^P_{\tv} \simeq \T_{\tv_{EGBP}}$ and, since the homeomorphism $T
\arrow
\td{T}$ takes proximal pairs to proximal pairs, we have

\begin{prop} \label{chain}  Suppose that $\vp$ is strong Pisot and satisfies
GCC and the no  cycle condition.  Then  $$
\T_{\vp_{EBP}} \simeq
\T^P_{\vp}
\simeq  \T^P_{\tv}
\simeq
\T_{\tv_{EGBP}}.$$
\end{prop}

The definition of  {\em weak
equivalence} for substitutions appears in \cite{bd1}.  For substitutions
$\tv$ and $\tilde{\psi}$, this reads as follows:  $\tv$ and $\tilde{\psi}$
are {\em weak equivalent}, denoted by 
$\tv \sim_w \tilde{\psi}$, if for $i
\in
\N$, there are
$n_i, m_i
\in
\N$ and morphisms $\tau_i : \A_{\tv} \arrow (\A_{\td{\psi}})^*$, $\sigma_i :
\A_{\td{\psi}} \arrow (\A_{\td{\vp}})^*$ so that $\tv^{m_i} = \sigma_i
\tau_i$ and
$\td{\psi}^{n_i} = \tau_i \sigma_{i+1}$.

\begin{picture}(25, 80)(-130, -40)

\put(-35, 0){$\sigma_1$}
\put(-10, 0){$\tau_1$}
\put(35, 0){$\sigma_2$}
\put(76, 0){$\tau_2$}
\put(110, 0){$\sigma_3$}
\put(-25, 20){$(\tilde{\vp})^{m_1}$}
\put(45, 20){$(\tilde{\vp})^{m_2}$}
\put(16, -32){$(\tilde{\psi})^{n_1}$}
\put(84, -32){$(\tilde{\psi})^{n_2}$}
\put(130, 0){$\ldots$}

\put(-20,-15){\vector(-1,1){30}}
\put(15,15){\vector(-1,-1){30}}
\put(50,-15){\vector(-1,1){30}}
\put(85,15){\vector(-1,-1){30}}
\put(120,-15){\vector(-1,1){30}}
\put(5,15){\vector(-1,0){45}}
\put(40,-15){\vector(-1,0){45}}
\put(110,-15){\vector(-1,0){45}}
\put(75,15){\vector(-1,0){45}}

\put(175,0){(2)}

\end{picture}

The next result follows immediately from known results. However, we make use
of the details of the argument  in the proof of Theorem  \ref{homess} and
include them  for completeness.

\begin{lemma} \label{weak} Suppose that $\vp$ and $\psi$ are strong Pisot and
satisfy both GCC and the no cycle condition.  Let 
$h:
\T_{\vp}
\arrow
\T_{\psi}$ be a homeomorphism.  There are one-cut rewritings $\tv$ of $\vp$
and
$\td{\psi}$ of
$\psi$ so that $\tv \sim_w \tilde{\psi}$ and  $$\inv f_{\td{\vp}}
\simeq
\T_{\td{\vp}} \simeq \T_\vp \simeq \T_\psi \simeq \T_{\td{\psi}} \simeq \inv
f_{\td{\psi}}.
$$ 

\end{lemma}

Proof: Let $h: \T_\vp \arrow \T_\psi$ be a homeomorphism.  Assume that
$h$ is orientation preserving; otherwise replace $h$ by its {\em reverse} (see
\cite{bd1}).    Also, we may assume that
$\vp$ and
$\psi$ are such that all asymptotic composants in $\T_\vp$ and $\T_\psi$ are fixed
by inflation and substitution (\cite{bd1}).  Since  $h$ takes
asymptotic composants  to asymptotic composants, and any composant fixed
under inflation and substitution contains a tiling  fixed under inflation and
substitution, we
 may modify $h$ by an isotopy to a homeomorphism $h'$ so that (1) for $T \in
\T_\vp$,  $h'(T) = h(T) + t$, where
$t = t(T)$, and (2) for some $T \in \T_\vp$, ${\Phi}(T) = T$ and
${\Psi}(h'(T)) = h'(T)$.  To simplify notation, we assume that $h$ itself has
this property.

We use these fixed tilings to determine the one-cut rewritings $\tv$ of $\vp$
and $\td{\psi}$ of
$\psi$. The roses
$R_{\tilde{\vp}}$ and
$R_{\tilde{\psi}}$  are formed by taking the disjoint unions of the
collections of prototiles and identifying all of the endpoints to a single
branch point $b$.  The prototiles $P_a$, for $a \in
\A_{\tilde{\vp}}$ or $a \in
\A_{\tilde{\psi}}$, then become the petals of the roses. The rose maps
$f_{\tilde{\vp}} : R_{\tilde{\vp}} \arrow R_{\tilde{\vp}}$ and
$f_{\tilde{\psi}} : R_{\tilde{\psi}} \arrow R_{\tilde{\psi}}$   fix
$b$ and map the petals following the patterns described by
$\tilde{\vp}$ and $\tilde{\psi}$, locally stretching arc length by a factor
of
$\lambda = \lambda_{\tilde{\vp}}$ or $\lambda = \lambda_{\tilde{\psi}}$.
Since $\tv$ and $\td{\psi}$ are proper, it follows from Theorem 4.3 of
\cite{anderson} that
$$\inv f_{\td{\vp}}
\simeq
\T_{\td{\vp}} \simeq \T_\vp \simeq \T_\psi \simeq \T_{\td{\psi}} \simeq \inv
f_{\td{\psi}} .$$ The homeomorphism $\hat{p}_{\tilde{\vp}} : \T_{\tilde{\vp}}
\arrow \inv f_{\tilde{\vp}}$ is defined by $$\hat{p}_{\tilde{\vp}} (T) =
(p_{\tilde{\vp}}(T), p_{\tilde{\vp}}(\td{\Phi}^{-1}(T)), 
p_{\tilde{\vp}}(\td{\Phi}^{-2}(T)),
\ldots ),$$ where $p_{\tilde{\vp}}(T) = s \in P_a \subseteq
R_{\tilde{\vp}}$ provided
$T_0$, the tile in
$T$ containing 0, is the translated prototile
$P_a$ with $T_0 = P_a - s$. It follows that  the homeomorphism of $\inv
f_{\tv}$ with $\inv f_{\td{\psi}}$ takes
$(b, b, \ldots) $ to $(b, b, \ldots)$, where $b$ denotes the branch point in
both $R_{\tv}$ and
$R_{\td{\psi}}$.  According to the proof of Theorem 1.16 of \cite{mjv}, $\tv$
is weakly equivalent with
$\td{\psi}$.

 \qed

\begin{theorem} \label{homess} Suppose that $\vp$ and $\psi$ are strong
Pisot and  satisfy both GCC and the no cycle condition.  Then
$\T_{\vp} \simeq
\T_{\psi}
 $ if and only if $ \T_{\vp_{EBP}} \simeq  \T_{\psi_{EBP}}  $.
\end{theorem}

Proof: Suppose that $ \T_{\vp_{EBP}} \simeq
\T_{\psi_{EBP}}  $. The original tiling spaces  $\T_{\vp} $ and $ \T_{\psi}
 $ sit in $ \T_{\vp_{EBP}} $ and $ \T_{\psi_{EBP}}  $ as distinguished
subspaces; for example, they are the unique subspaces irreducible with
respect to the property of being indecomposable.  It follows that  $\T_{\vp}
\simeq  \T_{\psi}
 $.

Now suppose that $\T_{\vp} \simeq  \T_{\psi}
 $.  According to Proposition~\ref{chain}, in order to show that $
\T_{\vp_{EBP}}
\simeq  \T_{\psi_{EBP}}  $, it is enough to show that
$\T_{\tv_{EGBP}} \simeq \T_{\td{\psi}_{EGBP}}$ for an appropriate choice of
$\tv,
\td{\psi}$.  We take  $\tv$ and $
\td{\psi}$ to be as in the proof of Lemma \ref{weak}, so that $\tv$ is weakly
equivalent with
$\td{\psi}$.

 Let $s_i$  and $t_i$ be the matrices (abelianizations) of the morphisms
$\sigma_i$ and $\tau_i$ (see (2)), so that
$(A_{\tv})^{m_i} = s_i t_i$ and
$(A_{\td{\psi}})^{n_i} = t_i s_{i+1}$.

Suppose that $\left( \begin{array}{c} \td{u}
\\ \td{v}
\end{array} \right)$ is a geometrically balanced pair for $\tv$, and let
$\omega_{L, \td{\vp}}$ denote a left Perron-Frobenius eigenvector for
$A_{\tv}$.  Then
$|\td{u}|_g = |\td{v}|_g$, so
$\Sigma|\td{u_i}|_g = \Sigma |\td{v_i}|_g$ and $\langle l(\td{u}), \omega_{L,
\td{\vp}}\rangle
  =
\langle l(\td{v}), \omega_{L, \td{\vp}}\rangle $ (where $\langle, \rangle$
denotes the usual Euclidean inner product).  That is, 
$\langle l(\td{u}) - l(\td{v}) , \omega_{L, \td{\vp}}\rangle = 0$,  which
implies that
$l(\td{u}) - l(\td{v})  \in E^s_{\tv}$ and  $(A_{\tv})^n (l(\td{u}) -
l(\td{v}))
\arrow 0$ as
$n \arrow
\infty$.  But
$l(\td{u}) - l(\td{v})  \in \Z^{\td{d}}$, so $(A_{\tv})^n (l(\td{u}) -
l(\td{v}))  = 0$ for some $n$, and $l(\td{u}) - l(\td{v}) \in N_{\tv}$, the
generalized null space of
$A_{\tv}$.

Choose $j$  large enough so that $ m := \Sigma_{i = 2}^j m_i \geq \td{d}$,
and define  $n = \Sigma_{i = 1}^{j-1} n_i$.  The fact that
$t_1 (A_{\tv})^m =  (A_{\td{\psi}})^n t_j$ implies that $t_j ( N_{\tv})
\subseteq N_{\td{\psi}}$.  Therefore $t_j (l(\td{u})) - t_j (l(\td{v})) \in
N_{\td{\psi}}$ for sufficiently large
$j$.  That is, $\langle t_j (l(\td{u})) - t_j (l(\td{v})), \omega_{L,
\td{\psi}}\rangle = 0$, so that
$|\tau_j(\td{u})|_g = |\tau_j(\td{v})|_g $.  In other words, for sufficiently
large $j$, $\tau_j$ takes geometrically balanced pairs for $\tv$ to
geometrically balanced pairs for
$\td{\psi}$.

Similarly,  for large enough $j$, $s_j (N_{\td{\psi}}) \subseteq N_{\tv}$,
and  $\sigma_j$ takes geometrically balanced pairs for $\td{\psi}$ to
geometrically balanced pairs for
$\td{\vp}$. Thus we have a weak equivalence between  $ \td{\vp}_{GBP}$ and $
\td{\psi}_{GBP} $.  It is clear that this restricts to a weak equivalence
between  $
\td{\vp}_{EGBP}$ and $
\td{\psi}_{EGBP} $.  A weak equivalence between substitutions induces a
homeomorphism between their tiling spaces (\cite{bd1}) so that  $
\T_{\td{\vp}_{EGBP}} \simeq
\T_{\td{\psi}_{EGBP}}  $.  \qed

\begin{theorem} \label{homprox} Suppose that $\vp$ and $\psi$ are strong
Pisot substitutions that satisfy both GCC and the no cycle condition, and let
$h: \T_\vp \arrow \T_\psi$ be a homeomorphism.  If $T, T'$ are proximal in
$\T_\vp$, there is  $t_0$ so that
$h(T)$ and $h(T') + t_0$ are proximal in
$\T_\psi$.
\end{theorem}

Proof:  Suppose that $h: \T_\vp \arrow \T_\psi$ is a homeomorphism, and let
$\tilde{\vp}$ and
$\tilde{\psi}$ be  one-cut rewritings of
$\vp$ and $\psi$ associated with $h$ defined as in the  proof
of Lemma
\ref{weak}.  We also use
$h:
\T_{\tilde{\vp}} \arrow
\T_{\tilde{\psi}}$ to denote the homeomorphism between $\T_{\tilde{\vp}} $
and $
\T_{\tilde{\psi}}$ induced by the homeomorphisms $T \arrow \td{T}$ associated
with each of
$\vp$ and $\psi$. Other notation in the following argument will serve to
avoid confusion.

  Recall that the homeomorphism $\hat{p}_{\tilde{\vp}} : \T_{\tilde{\vp}}
\arrow \inv f_{\tilde{\vp}}$ is defined by $$\hat{p}_{\tilde{\vp}} (T) =
(p_{\tilde{\vp}}(T), p_{\tilde{\vp}}(\td{\Phi}^{-1}(T)), 
p_{\tilde{\vp}}(\td{\Phi}^{-2}(T)),
\ldots ),$$ where $p_{\tilde{\vp}}(T) = s \in P_a \subseteq
R_{\tilde{\vp}}$ if
$T_0$, the tile in
$T$ containing 0, is the translated prototile
$P_a$ with $T_0 = P_a - s$.  Given $T = \{T_i\}_{i = - \infty}^\infty$, let
$\underline{w}(T) = \ldots w_{-1}w_0w_1 \ldots$ denote the bi-infinite word
representing $T$ (that is,   $w_i = j$ if and only if $T_i$ is of type
$j$).  Recall that the shift class of $\underline{w}(T)$,  $[
\underline{w}(T)]$, is the  pattern of the composant of
$T$.  Define the $k^{th}$ projection map $\pi_k : \inv f_{\tilde{\vp}}
\arrow R_{\tilde{\vp}}$ by
$\pi_k(x_0, x_1, \ldots) = x_k$.  Then the path $t \mapsto \pi_k(\hat{p}(T +
t))$ winds around
$ R_{\tilde{\vp}} $ following the pattern of the composant of
$\td{\Phi}^{-k}(T)$.

The weak equivalence between  $\tilde{\vp}$ and
$\tilde{\psi}$ induced by the homeomorphism $\hat{h} :=
\hat{p}_{\tilde{\psi}} \circ \ h \circ
\hat{p}^{-1}_{\tilde{\vp}}$ arises as follows (see \cite{mjv} for details).
For $a
\in
\A_{\tilde{\vp}}$, let
$\overset\circ{P_a} = P_a \setminus \{b\} \subseteq R_{\tilde{\vp}}$.   For
each
$k \in \N$ and $a \in
\A_{\tilde{\vp}}$, the set
$\overline{\pi^{-1}_k (\overset\circ{P_a}) }$
 is homeomorphic with the product of a Cantor set ($\pi^{-1}_k (\{x\})$,
where $x \in \overset\circ{P_a}$) and an arc. The larger is $k$, the longer
and skinnier is this product and the closer are all of its endpoints
$\overline{\pi^{-1}_k (\overset\circ{P_a})} \cap \pi^{-1}_k (\{b\})$ to the
point $\underline{b} = (b, b, \ldots )$.  A similar statement can be made for
$\inv f_{\tilde{\psi}}$.

Thus, there is an $l_1$ large enough so that, for each $a \in
\A_{\td{\psi}}$, the points $\pi_0(\hat{h}^{-1}(\overline{\pi^{-1}_{l_1}
(\overset\circ{P_a})} \cap \pi^{-1}_{l_1} (\{b\})))$ in $R_{\tilde{\vp}}$ are
all close to $b$, and all the arcs of $\hat{h}^{-1}(\overline{\pi^{-1}_{l_1}
(\overset\circ{P_a})})$ map under $\pi_0$ around $R_{\tilde{\vp}}$ in the
same pattern, which we denote by
$\sigma_1(a) \in \A^*_{\tv}$.

There is now an $m_1$ large enough so that  for each $a \in
\A_{\tv}$,  the points
$\pi_{l_1}(\hat{h}(\overline{\pi^{-1}_{m_1} (\overset\circ{P_a})} \cap
\pi^{-1}_{m_1} (\{b\})))$ are all close to $b$ in $R_{\tilde{\psi}}$, and the
arcs of $\hat{h}(\overline{\pi^{-1}_{m_1} (\overset\circ{P_a})}) $ all map in
the same well-defined pattern, which we label
$\tau_1(a) \in \A^*_{\td{\psi}}$, around $R_{\tilde{\psi}}$, etc.  We see
from this description that if $[\underline{w}]$ is the pattern of the
composant of
$\td{\Phi}^{-m_1}(T)$, then $[\tau_1(\underline{w})]$ is the pattern of the
composant of $\td{\Psi}^{-l_1} (h(T))$.

Suppose that $T$ and $T'$ are proximal in $\T_\vp$.  Then $\td{T}$ and
$\td{T}'$ are proximal in $\T_{\tv}$, and hence so are
$\td{\Phi}^{-m_1}(\td{T})$ and
$\td{\Phi}^{-m_1}(\td{T}')$.  Thus the pair of words  $
\left(
\begin{array}{c} \td{\Phi}^{-m_1}(\td{T}) \\ \td{\Phi}^{-m_1}(\td{T}')
\end{array}
\right)$  factors as a bi-infinite product of essential geometrically
balanced pairs for $\tv$.  Apply $\tau_1$ to this product; the result factors as a
bi-infinite product of essential geometrically balanced pairs for
$\td{\psi}$ (see the proof of Theorem \ref{homess}).  Thus the patterns
$[\tau_1 \underline{w}(\td{\Phi}^{-m_1}(\td{T}))]$ and
$[\tau_1 \underline{w}(\td{\Phi}^{-m_1}(\td{T}'))]$, appropriately shifted,
balance geometrically.  This means that $\td{\Psi}^{-l_1} (h(T))$ and
$\td{\Psi}^{-l_1} (h(T')) + t_1$ are proximal in $\T_{\td{\psi}}$ for some
$t_1$, and hence so are
$h(T)$ and $h(T') + t_0$, where $t_0 = (\lambda_{\td{\psi}})^{l_1}  t_1$.
$\qed$

\

\begin{example} \label{diffprox} Using proximality to distinguish tiling
spaces
\end{example}

Define $\vp$ and $\psi$ as follows:
$$ \vp(1) = 1112211122111221212, \vp(2) = 1112212; $$
$$  \psi(1) = 1112211121212121212, \psi(2) = 1112212. $$

The substitution  $\vp$ is the second iterate of that
considered in Examples \ref{esspairs}, \ref{rewriting} and \ref{proximal}, and the
substitution
$\psi$ is obtained by modifying $\vp$ slightly. 

The basic
 structure of the asymptotic composants is the
same: A proper substitution on two symbols can have
at most four asymptotic composants (\cite{charles}), and 
$\T_\vp$ and $\T_\psi$ both have a pair of backward asymptotic and a pair of
forward asymptotic composants.   For  
$\vp$, these are represented by the bi-infinite words (spaced to indicate
balanced pairs):
$$ \begin{array}{c}  \ldots  1\ 2\ 1\ 1\ 1\ 2\ 2\  1\ \underline{1} 12\  2 1\ 
112\ 2 1\  2 
\ldots\\ \ldots  1\ 2\ 1\ 1\  1 \ 2\  2\  1\ 2  .  1 1\  1 2\  2 11\ 12 \ 2
\ldots \end{array} $$ and $$ \begin{array}{c}  \ldots 2\ 2 11\ 12\ 2\ 1\ 112\ 2
\underline{1}\ 2\ 1\ 2 \ 1\ 1\ 1\ 
\ldots\\ \ldots 1\ 112\ 21\ 2\ 1\ 2  1 1\ 1
\underline{2}\ 2\ 1\ 2 \ 1\ 1\ 1\ 
\ldots \end{array}   $$
and for $\psi$ by 
$$  \begin{array}{c} \ldots 1\ 1\ 1\ 2\ 2\ 1\ \underline{1}12\ 1\ 2\ 12\ 1\ 21\ 2
\ldots
\\
\ldots 1\ 1\ 1\ 2\ 2\ 1\ 2 . 11\ 1\ 2\ 21\ 1\ 12 \ 2 \ldots 
\end{array} $$ and $$  \begin{array}{c} \ldots 1\ 12\ 12 \ 1\ 2\underline{1}\ 2\
1\ 2
\ 1\ 1\ 1\ 
\ldots \\
\ldots 1\ 21\ 2 1\ 1\ 1\underline{2}\ 2\ 1\ 2 \ 1\ 1\ 1\  \ldots 
\end{array} .$$

It is easy to check that $\vp$ and
$\psi$ are strong Pisot.  All strong Pisot substitutions on two symbols  satisfy
GCC (see
\cite{bd2}, \cite{holsol}, and
\cite{mjarek} for the unimodular case).

Tilings, asymptotic composants,  and essential balanced pairs are identical for 
all iterates of a substitution, and we saw in Examples \ref{esspairs} and 
\ref{proximal} that $\vp$ has at least three essential balanced pairs
which  generate the backward asymptotic (and forward proximal) tilings, $T$ and
$T'$, and a third tiling $T''$ proximal in both
directions with each of  the first two.  This can be seen from $\vp_{EBP}$, which 
includes at least the information (recall that $\vp$ is  the second iterate of
the substitution in Example \ref{esspairs}):
$$\vp_{EBP}(a) = 11122 1 b \overline{a} b
\overline{a} 2 1\overline{b}a2 12,$$ $$ \vp_{EBP}(b) = 11122 1 b \overline{a} b
\overline{a} 2 1\overline{b}a2 1 b \overline{a} b
\overline{a} 2 1\overline{b}a2 12,$$ $$ \vp_{EBP}(c) = 11122 1 b \overline{a} b
\overline{a} 2 1\overline{b}a2 1b \overline{a} c 11122
1\overline{b} a2 12,$$ along with the definition of $\vp_{EBP}$ on duals.  In
particular, $\vp_{EBP}$ has two additional backward  asymptotic composants
indicating that
$T''$ is proximal to $T$ and $T'$ in $\T_\vp$:   
$$  \begin{array}{c} \ldots 1 b \overline{a} \underline{b}
\overline{a} 2 1
\ldots \\
\ldots 1 b \overline{a} \underline{c} 1
\overline{a} 2 \ldots 
\end{array} $$  along with the three that capture the two original backward
asymptotic composants:
$$ \begin{array}{c}  \ldots  111221 \underline{1} 12 2 
\ldots\\ \ldots  1112212  .  122
\ldots \\ \ldots 111221\underline{b} \overline{a}b\overline{a} \ldots  .
\end{array} $$

Suppose that $ \T_\vp \simeq
 \T_\psi$.   According to Theorem \ref{homess},  $
\T_{\vp_{EBP}} \simeq  \T_{\psi_{EBP}}  $. As in the proof of Lemma \ref{weak}, we
assume the homeomorphism 
 is  orientation preserving. Since such a homeomorphism  must take backward
asymptotic composants to backward asymptotic composants, $ \T_{\psi_{EBP}}  $
must have at least five backward asymptotic composants.

It follows from \begin{lemma} \label{onlypairs} The  balanced pairs $\larr 1 \\ 1
\rarr$, 
$\larr 2 \\ 2
\rarr$,  $\larr 21 \\ 12
\rarr$, 
$\larr 211
\\ 122
\rarr$, and their duals are the only essential balanced pairs for $\psi$.
\end{lemma} 

{\parindent=0pt  that}
 the substitution $\psi_{EBP}$  is entirely  given by $$
\begin{array}{c} 
  \psi_{EBP}(a)
= 1 1 1 2 2 1 b 1 2 a 1 \overline{a} \overline{a} 1 a 2 1 2 \\ \psi_{EBP}(b) = 1 1
1 2 2 1 b 1 2 a 1 \overline{a} \overline{a} \overline{a} 1 a 2 1 b 1 2 a 1
\overline{a}\overline{a} \overline{a} 1 a 2 1 2 \end{array}$$ and the implied
definition on trivial balanced pairs and the duals of $a, b$.  The only backward
asymptotic words for
$\psi_{EBP}$ are $$  \begin{array}{c} \ldots 111221 \underline{1}121 
\ldots
\\
\ldots 111221 \underline{b} 1 2 a  \ldots \\ \ldots 111221\underline{2} 11 1
\ldots
\end{array}  $$ which code  the original backward asymptotic
composants.  That is, $\T_{\psi_{EBP}}$ has only three backward asymptotic
composants, and
$\T_\vp
\nsim
\T_\psi$.   

The proof of Lemma \ref{onlypairs}  appears in the Appendix.
\qed 

\

Suppose that $\vp$ and $\psi$ are strong Pisot and satisfy both GCC and the
no cycle condition. If $\T_\vp \simeq \T_\psi$, then 
$\T_\vp^P$ is homeomorphic to $\T_\psi^P $ under a homeomorphism that maps
$\T_\vp$ to $\T_\psi$, hence 
$\T_\vp^P /
\T_\vp \simeq \T_\psi^P / \T_\psi$ (we are identifying $\T_\vp$ and $\T_\psi$ 
with the `diagonals' in $\T_\vp^P$ and $\T_\psi^P$).  The space
$\T_\vp^P /
\T_\vp$ has a local product structure everywhere but at $[\T_\vp]$.  The
element
$[\T_\vp]$ itself is the center of an $m$-od, where $m = 2(\#$ asymptotic
pairs). Let $R_{\vp_{EBP}}$ be the rose associated with $\T_\vp^P \simeq
\T_{\vp_{EBP}}$, and let $f_{\vp_{EBP}}: R_{\vp_{EBP}} \arrow R_{\vp_{EBP}}$
be the rose map. Collapsing $ R_{\vp}$ (which is embedded in $  
R_{\vp_{EBP}}$) to the branch point induces a substitution $\vp_P$ on
just the symbols in $\A_{EBP}$ that correspond to nontrivial essential
balanced pairs and the  map
$f_{\vp_P} :  R_{\vp_P} := R_{\vp_{EBP}}/ R_\vp
\arrow R_{\vp_P}  $.  Furthermore,  $\T_\vp^P/\T_\vp \simeq \inv f_{\vp_P}$.

Assuming still that $\vp$ and $\psi$ are strong Pisot and satisfy both GCC
and the no cycle condition, suppose that
$\T_\vp \simeq \T_\psi$.  Any homeomorphism of $\T^P_\vp$ with
$\T^P_\psi$ not only takes $\T_\vp $ to $\T_\psi$ but must also take arc
components of
$\T^P_\vp$ that are asymptotic to asymptotic composants of $\T_\vp$ to arc
components of $\T^P_\psi$ that are asymptotic to asymptotic composants of
$\T_\psi$ (again,  $\T_\vp$ and $\T_\psi$ are identified
with the `diagonals' in $\T_\vp^P$ and $\T_\psi^P$).  Let $\T^A_\vp$ be the minimal subcontinuum of $\T^P_\vp$ that
contains all arc components of $\T^P_\vp$ that are asymptotic to asymptotic
composants of $\T_\vp$.  Thus $\T_\vp \simeq \T_\psi$ implies
 $\T^A_\vp \simeq \T^A_\psi$.

Note that $\larr T'\\ T''\rarr \in \T^P_\vp$ is asymptotic to $\larr T\\
T \rarr \in \T_\vp \subset \T^P_\vp$ if and only if $T, T', T''$ are all
asymptotic (in the same direction).  Thus the pairs  $\larr T'\\ T''\rarr
\in \T^A_\vp$ are precisely those proximal pairs all of whose bubbles come
from asymptotic pairs (see Proposition \ref{bubble}). More precisely, let
\begin{eqnarray*} ABP_\vp &= &  \{ \larr u \\ v \rarr: \larr u \\ v
\rarr \textrm{ is a trivial balanced pair for } \vp \textrm{ or } \larr u
\\ v \rarr \\ & &
\textrm{ is obtained from a bubble in an asymptotic pair for } \vp \}.
\end{eqnarray*} 

As $\Phi$ takes asymptotic pairs to asymptotic pairs, for $\larr u \\ v
\rarr \in ABP_\vp$,  
$\larr
\vp(u)
\\
\vp(v) \rarr$ can be factored as a product of elements of $ABP_\vp$; this
defines $\vp_{ABP}: ABP_\vp \arrow (ABP_\vp)^*$.  We see that
$\T_{\vp_{ABP}} \simeq \T^A_\vp$.  Letting $\vp_A$ be the substitution on
the nontrivial elements of $ABP_\vp$ that is the composition of $\vp_{ABP}$
with the morphism that forgets the trivial balanced pairs, and letting  $f_{\vp_A}
: R_{\vp_A} \arrow R_{\vp_A}$ be the associated rose map, we have (just as in the
preceding paragraph) that $\T^A_\vp/\T_\vp \simeq \inv f_{\vp_A}$.  We have
almost completed the proof of 
\begin{prop} \label{reducedth} Suppose that $\vp$ and $\psi$ are strong Pisot
and satisfy both GCC and the no cycle condition.  If $\T_\vp \simeq
\T_\psi$, then:
\\(i)  $\inv f_{\vp_P} \simeq \inv f_{\psi_P}$ and \\ (ii)  $\inv f_{\vp_A}
\simeq \inv f_{\psi_A}$.
\\ Moreover (i) is equivalent to $ \vp_P \sim_w
\psi_P$ and (ii) is equivalent to $ \vp_A \sim_w
\psi_A$.  (The definition of $\sim_w$ precedes Lemma \ref{weak}.)

\end{prop}

Proof:  (i) and (ii) follow from the discussion preceding this
proposition.  The branch points of 
$\inv f_{\vp_P}$,
$\inv f_{\psi_P}$, $\inv f_{\vp_A}
$ and $\inv f_{\psi_A}$ are distinguished, so by \cite{mjv}, $\inv f_{\vp_P}
\simeq \inv f_{\psi_P}$ if and only if $ \vp_P \sim_w
\psi_P$ and $\inv f_{\vp_A}
\simeq \inv f_{\psi_A}$ if and only if $ \vp_A \sim_w
\psi_A$.  \qed

\

{\parindent=0pt{\em Remark}:}  It can be tedious to identify all essential
balanced pairs, even for relatively simple substitutions.  At the same time,
the alphabet for the substitution
$\vp_A$ on nontrivial balanced  pairs arising from bubbles formed by
asymptotic pairs is no more difficult to determine than the asymptotic pairs
themselves.  As a result, the use of  (ii) of Proposition
\ref{reducedth}  to distinguish nonhomeomorphic tiling spaces $\T_\vp$
and $\T_\psi$ is more straightforward when it applies.  We illustrate this
with the next example.

\begin{example} Distinguishing tiling spaces using the reduced substitution
on balanced pairs
\end{example}

As part of his program to classify hyperbolic one-dimensional attractors (up
to topological conjugacy), R.F.  Williams (\cite{wil1}) sought to determine
the shift equivalence classes of all (there are 46) substitutions on two
letters that are proper and whose abelianizations have characteristic
polynomial $x^2 - 3x - 2$.  Shift equivalence of proper  substitutions is
equivalent to topological conjugacy of the corresponding inflation and
substitution homeomorphisms of the tiling spaces.  Williams reduced the
problem to the consideration of four particular substitutions, two of which:
$$\vp(1) = 11221,
\vp(2) = 1,$$ and   $$\psi(1) = 112222, \vp(2) = 12,$$ were finally shown not
to be shift equivalent in \cite{da}.  We show here that 
$\T_\vp $ and $\T_\psi$ are not even homeomorphic.    This completes the
topological classification of the spaces arising from the characteristic
polynomial $x^2 - 3x - 2$:  two of these spaces are homeomorphic if and
only if their inflation and substitution homeomorphisms are conjugate, and
there are exactly three topological equivalence classes (see
\cite{mrichard}).

Each of $\T_\vp $ and $\T_\psi$  has one pair of forward and one pair
of backward asymptotic composants.  Since a tiling space for a proper substitution
on two letters has at most four asymptotic composants (\cite{charles}), there are
no others.

A proper substitution has no cycles of periodic words.  Also, since $\vp$ and
$\psi$ are strong Pisot substitutions on two letters, they satisfy GCC
(see \cite{bd2}, \cite{holsol}, and \cite{mjarek} for the unimodular case--the
nonunimodular case is a straightforward generalization).  It is tedious to ensure
that  all essential balanced pairs have been identified for either of
$\vp$ and
$\psi$, so we make use of Proposition \ref{reducedth} and work only with 
essential balanced pairs associated with asymptotic composants.

 For
$\vp$, the only essential balanced pair associated with asymptotic
composants is $a:= \larr 1 1 2 2 \\2 2 11\rarr$
 and its dual $\overline{a} $.   If $1 := \larr 1 \\ 1\rarr$ and $2 := \larr 2 \\
2\rarr$, then
$\vp_{ABP}(i) = \vp(i)$ for $i = 1, 2$, and $$ \vp_{ABP}(a) = 1 1 \overline{a}
1\overline{a} 1, 
\vp_{ABP}(\overline{a} ) = 11a1a1.
$$  The reduced substitution is $$\vp_A(a) = \overline{a} \overline{a} ,
\vp_A(\overline{a} ) = aa,$$ and the inverse limit of the rose map, $\inv
f_{\vp_A}$, is homeomorphic with a pair of dyadic solenoids joined at a point.

As for $\psi$, there are ten essential balanced pairs associated with asymptotic
composants:  $a:= \larr 12 \\ 21\rarr$, $b:= \larr 2221 \\ 1222\rarr$, $c:= \larr
1211222\\ 2212121\rarr$, $d:= \larr 22121\\ 11222\rarr$, 
$e:=
\larr 1 1 2 2
\\2 2 11\rarr$, and their duals.  A computation yields: $$\psi_A(a) = ab,
\psi_A(b) = \overline{a}c, \psi_A(ac) = ada\overline{e}\overline{c}, \psi_A(d) =
\overline{a}ec, 
\psi_A(e) =
a\overline{e}d$$ (and the dual statements).  A quick check shows that $\psi_A$ is
primitive and aperiodic, so $\inv f_{\psi_A}$ is an indecomposable
continuum, thus $\inv f_{\vp_A} \nsim \inv f_{\psi_A}$.  (Alternatively,
$\check{H}^1(\inv f_{\vp_A}) \simeq \dir \larr 0 \ 2 \\ 2 \ 0 \rarr \simeq \Z
\oplus \Z$, while $\check{H}^1(\inv f_{\psi_A})$ is homeomorphic with the direct
sum of ten copies of $\Z$, as $A$ is nonsingular.)  Thus, by Proposition
\ref{reducedth} (ii), $ \T_\vp \nsim \T_\psi  $. \qed

\section{GCC if and only if proximality is closed}\label{sec-closed}
The main result of this section is  Theorem \ref{closedGCC}, in which we show
that a strong Pisot substitution  $\vp$ satisfies GCC if and only if proximality
is a closed relation on $\T_\vp$.
\begin{lemma} \label{dbyd} Suppose that $A$ is a nonsingular, hyperbolic,
integer
$d
\times d$ matrix and $F_A: \T_d \arrow \T_d$ is the associated Anosov
endomorphism.  Suppose also that there is a  compact metric space $X$
with homeomorphism $F: X \arrow X$ and a covering map
$c: X
\arrow
\inv F_A$ that semi-conjugates $F$ with the shift homeomorphism $\hat{F}_A :
\inv F_A \arrow \inv F_A$.  Then there are $d \times d$ integer matrices $B$ and
$R$ with
$RB = AR$ and a homeomorphism $h: X \arrow \inv F_B$ so that the diagram 

{\parindent=0pt
\begin{picture}(25, 80)(-130, -40)

\put(80, -8){$\hat{F}_R$}
\put(23, -32){$\inv F_A$}
\put(-10, -8){$c$}

\put(-20, 20){$ X$}
\put(70, 20){$\inv F_B $}
\put(-10,10){\vector(1,-1){30}}
\put(-0,25){\vector(1,0){65}}
\put(35, 30){$h$}

\put(85,10){\vector(-1,-1){30}}

\end{picture}  

commutes,} 
 where $\hat{F}_R$ is induced by $F_R: \cT_d \arrow \cT_d$.
\end{lemma}

Proof:  Let $c$ be an $m$-to-one covering map and let $r = \deg F_A = |\det
A|$.  Let $\pi_k: \inv F_A \arrow \cT_d$ be projection onto the $k^{th}$
coordinate.  Given $\delta > 0$, let $\{U_i\}$ be a finite cover of
$\cT_d$ by open $\delta/4$ balls with the property that if $U_i \cap U_l =
\emptyset$, then $\overline{U}_i \cap \overline{U}_l =
\emptyset$.  For each $k \in \N$, $F^{-k}_A(U_i) = U'_i \cup \ldots \cup
U_i^{r^k}$ is a disjoint union of topological balls.  For large enough $k$,
$diam ( \pi^{-1}_k(U^j_i)) < \delta$ for all $i, j$.  Then, for sufficiently
small $\delta$, $c^{-1}(\pi_k^{-1}(U_i^j)) = W_i^{j, 1} \cup \ldots \cup
W_i^{j, m} $ is a disjoint union with $c|W_i^{j, s} : W_i^{j, s} \arrow
\pi_k^{-1}(U_i^j)$ a homeomorphism for all $i, j, s$.  Define the relation
$\sim$ on $X$ by $x \sim y$ if and only if $x, y \in W_i^{j, s}$ for some $i,
j, s$ and $\pi_k \circ c(x) = \pi_k \circ c(y)$.  Note that if $x \in
W_i^{j, s} \cap W_l^{t, q}$, $y \in W_i^{j, s}$ and $\pi_k \circ c(x) =
\pi_k \circ c(y)$, then $ y \in  W_l^{t, q}$.  It follows that $\sim$ is an
equivalence relation.  Let $X_1 = (X/\sim)$, and define $p_1 : X_1 \arrow
\cT^d$ by $p_1([x]) = \pi_k(c(x))$.  Then $p_1$ is exactly $m$-to-one
everywhere, and if $p_1([x]) = p_1([y])$ with $[x] \neq [y]$, then there are
$s, i, j, q, t, l$ with $x \in W_i^{j, s} $, $y \in W_l^{t, q}$, and
$W_i^{j, s} \cap W_l^{t, q} = \emptyset$.  Since $\overline{W}_i^{j, s} \cap
\overline{W}_l^{t, q} = \emptyset$ by the assumption on $U_i$ and $U_l$,
the distance between $[x] $ and $[y]$ must be at least as large as the
minimum of all the minimum distances between pairs of disjoint compact sets
$\overline{W}_b^{c, a} $ and $\overline{W}_e^{f, g} $. That is, there is
$\eta > 0$ so that if $p_1([x] ) = p_1([y])$ and $[x] \neq [y]$, then
$d([x], [y]) > \eta$.  Thus $p_1$ is a covering map and since $X_1$ is
compact, $X_1 \simeq \cT^d$.

Now if $x \sim y$, say $x, y \in W_i^{j, s} $ with $\pi_k \circ c(x) =
\pi_k \circ c(y)$, then $$\pi_k \circ c(F(x)) = \pi_k \circ
\hat{F}_A(c(x)) = F_A(\pi_k(c(x))) $$ $$ \hspace*{1in}= F_A(\pi_k(c(y))) = 
\pi_k
\circ
\hat{F}_A(c(y)) = \pi_k \circ c(F(y)) .$$  We show that $F(x), F(y) \in
W_l^{t, q}$ for some $q, t, l$, so that $F(x) \sim F(y)$.  In order to
guarantee this, we will adjust $\delta$ (and $k$ correspondingly):  
 Let $\rho = \inf\{d(u, v) : u \neq v, c(u) =
c(v)\} > 0$, and choose $\epsilon > 0$  small enough so that if $d(x, y) <
\epsilon$, then
$d(F(x), F(y)) < \rho/2$. Now, for
$\delta > 0$ sufficiently small, $diam (W_i^{j, s}) < \epsilon$ for all $i,
j, s$. With this $\delta$ and $k$, let $\pi_k \circ
c(F(x))
\in U_l^t$.  Then $\pi_k \circ c(F(x)) = \pi_k \circ c(F(y)) \in U^t_l$,
hence for some $p, q$, $F(x) \in W_l^{t, p}$ and $F(y) \in W_l^{t, q}$. 
Since $x, y \in W_i^{j, s}$, $d(x, y) < \epsilon$.  It follows that $d(F(x), F(y)) <
\rho/2$, hence  $p = q$ and $F(x) \sim F(y)$.

Thus there is an induced map $F_1 : X_1 \arrow X_1$, and $$p_1 \circ
F_1([x]) = p_1 ([F(x)]) = \pi_k \circ c(F(x))  = \pi_k \circ
\hat{F}_A(c(x)) $$ $$= F_A \circ \pi_k(c(x))  = F_A \circ p_1([x]).$$ That
is, we have a commuting diagram of toral endomorphisms 

{\parindent=0pt
\begin{picture}(25, 95)(-140, -50)

\put(0, -7){$p_1$}
\put(75, -7){$p_1$}
\put(40, 25){$F_1$}
\put(40, -44){$F_A$}

\put(-3, 17){$X_1$}
\put(75, 17){$X_1$}
\put(-5, -33){$\cT^d$}
\put(75, -33){$\cT^d$}


\put(15,12){\vector(0, -1){37}}
\put(65,-30){\vector(-1,0){45}}
\put(70,12){\vector(0,-1){37}}
\put(65,20){\vector(-1,0){45}}

\end{picture}

It } now follows that there are integer matrices $R$ and $B$ (with $|\det R
| = m$) with $RB = AR$ and a homeomorphism $h_1 : X_1 \arrow \cT^d$ so that 
$F_B = h_1 \circ F_1 \circ h^{-1}_1$, $F_R \circ h_1 = p_1$ and the diagram 

{\parindent=0pt
\begin{picture}(25, 95)(-140, -50)

\put(-3, -7){$F_R$}
\put(75, -7){$F_R$}
\put(40, 25){$F_B$}
\put(40, -44){$F_A$}

\put(-3, 17){$\cT^d$}
\put(75, 17){$\cT^d$}
\put(-5, -33){$\cT^d$}
\put(75, -33){$\cT^d$}


\put(15,12){\vector(0, -1){37}}
\put(65,-30){\vector(-1,0){45}}
\put(70,12){\vector(0,-1){37}}
\put(65,20){\vector(-1,0){45}}

\end{picture}

commutes.}

Claim:  $\inv F_1 \simeq X$.

The commuting diagram

{\parindent=0pt
\begin{picture}(25, 140)(-140, -100)

\put(19, -20){$c$}
\put(92, -20){$c$}
\put(45, 25){$F$}
\put(40, -43){$\hat{F}_A$}

\put(13, 17){$X$}
\put(85, 17){$X$}
\put(-5, -50){$\inv F_A$}
\put(85, -50){$\inv F_A$}

\put(18,12){\line(0, -1){15}}
\put(18,-12){\vector(0, -1){25}}

\put(80,-47){\line(-1,0){18}}
\put(50,-47){\vector(-1,0){18}}
\put(90,12){\vector(0,-1){50}}
\put(80,20){\vector(-1,0){50}}


\put(8, 14){\vector(-1, -1){15}}
\put(80, 14){\vector(-1, -1){15}}
\put(8, -54){\vector(-1, -1){15}}
\put(80, -54){\vector(-1, -1){15}}

\put(-9, 5){$\pi$}
\put(62, 5){$\pi$}
\put(2, -67){$\pi_k$}
\put(74, -67){$\pi_k$}

\put(27, -4){$F_1$}
\put(15, -85){$F_A$}
\put(-28, -40){$p_1$}
\put(60, -30){$p_1$}

\put(-20, -12){$X_1$}
\put(52, -12){$X_1$}
\put(-20, -79){$\cT^d$}
\put(52, -79){$\cT^d$}

\put(-15,-17){\vector(0, -1){50}}
\put(47,-8){\vector(-1,0){50}}
\put(57,-17){\vector(0,-1){50}}
\put(47,-75){\vector(-1,0){50}}

\end{picture}

in } which $\pi: X \arrow (X/\sim) \  = X_1$ is the quotient map, induces the
commuting diagram   

{\parindent=0pt
\begin{picture}(25, 140)(-140, -100)

\put(19, -20){$\hat{c}$}
\put(92, -20){$\hat{c}$}
\put(45, 25){$F$}
\put(40, -44){$\hat{F}_A$}

\put(13, 17){$X$}
\put(85, 17){$X$}
\put(-5, -50){$\inv F_A$}
\put(85, -50){$\inv F_A$}

\put(18,12){\vector(0, -1){50}}
\put(80,-47){\vector(-1,0){50}}
\put(90,12){\vector(0,-1){50}}
\put(80,20){\vector(-1,0){50}}


\put(8, 14){\vector(-1, -1){15}}
\put(80, 14){\vector(-1, -1){15}}
\put(8, -54){\vector(-1, -1){15}}
\put(80, -54){\vector(-1, -1){15}}

\put(-9, 5){$\hat{\pi}$}
\put(62, 5){$\hat{\pi}$}
\put(2, -67){$\hat{\pi}_k$}
\put(74, -67){$\hat{\pi}_k$}

\put(27, -4){$\hat{F}_1$}
\put(15, -88){$\hat{F}_A$}
\put(-28, -40){$\hat{p}_1$}
\put(60, -30){$\hat{p}_1$}

\put(-35, -12){$\inv F_1$}
\put(52, -12){$\inv F_1$}
\put(-35, -79){$\inv F_A$}
\put(52, -79){$\inv F_A$}

\put(-15,-17){\vector(0, -1){50}}
\put(47,-8){\vector(-1,0){50}}
\put(57,-17){\vector(0,-1){50}}
\put(47,-75){\vector(-1,0){50}}

\end{picture}

Here} $\hat{\pi}_k = (\hat{F}_A)^{-k}$, $\hat{c}$ is exactly $m$-to-one, and
$\hat{p}_1$ is exactly $m$-to-one.  Thus $\hat{\pi}$ is a homeomorphism. 
Moreover, 

{\parindent=0pt
\begin{picture}(25, 200)(-120, -100)

\put(-3, -7){$\hat{h}_1$}
\put(75, -7){$\hat{h}_1$}
\put(40, 25){$\hat{F}_1$}
\put(40, -26){$\hat{F}_B$}

\put(-8, -52){$\hat{F}_R$}\put(80, -52){$\hat{F}_R$}
\put(-8, 40){$\hat{\pi}$}\put(85, 40){$\hat{\pi}$}

\put(-43, 0){$c$}\put(120, 0){$c$}
\put(40, 63){$F$}\put(40, -80){$\hat{F}_A$}

\put(-38, 55){$X$}
\put(-45, -70){$\inv F_A$}
\put(110, 55){$X$}
\put(110, -70){$\inv F_A$}

\put(-17, 17){$\inv F_1$}
\put(75, 17){$\inv F_1$}
\put(-17, -33){$\inv F_B$}
\put(75, -33){$\inv F_B$}

\put(15,12){\vector(0, -1){37}}
\put(20,-30){\vector(1,0){45}}
\put(70,12){\vector(0,-1){37}}
\put(20,20){\vector(1,0){45}}

\put(20,-30){\vector(1,0){45}}
\put(20,-30){\vector(1,0){45}}

\put(-33,45){\vector(0,-1){95}}
\put(117,45){\vector(0, -1){95}}
\put(-20,58){\vector(1, 0){125}}
\put(-5,-65){\vector(1, 0){110}}

\put(-20, 50){\vector(1, -1){15}}
\put(105, 50){\vector(-1, -1){15}}
\put(-7, -40){\vector(-1, -1){15}}
\put(90, -40){\vector(1, -1){15}}

\end{picture}

commutes.}  Letting $h = \hat{h}_1 \circ \hat{\pi}$, we have the desired
result.
\qed 

\

\begin{lemma} \label{closed} Suppose that $\vp$ is strong Pisot and that
$\sim_p$ is a closed relation on $\T_\vp \simeq \T_\vp^S$.  Then:
\\(1) $\gamma \in \T_\vp^S$ is forward proximal to $\gamma' \in \T_\vp^S$ if
and only if $\gamma$ is backwards proximal to $\gamma'$.
\\(2) If $\{B_i\} = \B$ is the collection of all irreducible balanced pairs
formed by proximal pairs in $\T_\vp^S$, then $\B$ is finite and each $B_i$
terminates with coincidence.
\\ (3)  $\sim_p$ is an equivalence relation.
\\(4)  If $\gamma \nsim_p \gamma' $ and $\hat{g}(\gamma) = \hat{g}(\gamma')$,
then $\gamma$ and $\gamma'$ do not share an edge.
\end{lemma}

Proof:  (1)  Suppose that $\gamma$ and $\gamma'$ are forward proximal, so
that there is a sequence $\{t_n\}$ with $d(\gamma - t_n \omega_R, \gamma'
- t_n \omega_R) \arrow 0$ as $t_n \arrow \infty$.  If $\gamma$ and $\gamma'$
are not backward proximal, there is $t_0 \in \R$ and $\epsilon > 0$ so that
$d(\gamma - t \omega_R, \gamma' - t \omega_R) \geq \epsilon$ for all $t \leq
t_0$.  Choose $\{s_n\}$ so that $s_n \arrow -\infty$ and $\gamma - s_n
\omega_R \arrow \eta$, $\gamma - s_n \omega_R \arrow \eta'$.  Since $\sim_p$
is closed, $\eta \sim_p \eta'$.   There is  $t \in \R$  so that $d(\eta - t
\omega_R, \eta' - t \omega_R) < \epsilon$.  Then, for sufficiently large
$n$, $d(\gamma - (s_n+t) \omega_R, \gamma' - (s_n + t) \omega_R) < \epsilon$
and $s_n + t \leq t_0$, a contradiction.  The converse can be proved by a
symmetric argument.

(2)  By (1), proximal pairs determine a bi-infinite product of irreducible
balanced pairs.  Suppose that there are arbitrarily long products of
nontrivial balanced pairs that arise in such factorizations.  Let $\gamma_i
\sim_p \gamma_i'$ with $\gamma_i$ and $\gamma_i'$ forming the products $W_i$
of nontrivial irreducible balanced pairs, where $W_i$ is centered on $E^s$
with $|W_i|_g \arrow \infty$.  Without loss of generality, $\gamma_i \arrow
\eta$ and $\gamma_i' \arrow \eta'$.  The relation $\sim_p$ is closed,  hence
$\eta \sim_p \eta'$.  But the definition of $\gamma_i$ and $\gamma_i'$
implies that $\eta \nsim_p \eta'$, a contradiction.  Thus $\B$ is finite. 
Since $\gamma \sim_p \gamma'$ implies that $\Phi(\gamma) \sim_p
\Phi(\gamma')$, each $B_i \in \B$ terminates with coincidence.

(3)  Consider the restriction of the substitution on balanced pairs,
$\vp_{BP}: BP(\vp) \arrow BP(\vp)$, to the set of irreducible (trivial and
nontrivial)   balanced pairs associated with proximal pairs.  According to
Lemma \ref{freq}, tiles of type corresponding to trivial balanced pairs have
frequency 1 in every element of $\T_{\tv}$. That is, coincidence in proximal
pairs of $\T_\vp^S$ occurs with frequency 1, hence proximality is transitive.

(4)  Suppose that $\hat{g}(\gamma) = \hat{g}(\gamma')$ and that $\gamma$ and
$\gamma'$ coincide along some segment $I$.  Let $m = cr_\vp$.  By Theorem
\ref{GCCeq}, for some $\underline{x} \in \inv F_A$,
$\hat{g}^{-1}(\underline{x})
$ consists of $m$ strands, any two of which are nowhere coincident.
  By
minimality of the flow on $\inv F_A$, in every preimage
$\hat{g}^{-1}(\underline{y})$, there are at least $m$ pairwise nowhere
coincident strands.  It follows that $$\G := \{\underline{y}:
\hat{g}^{-1}(\underline{y}) \textrm{ contains exactly } m \textrm{
pairwise nowhere coin. strands}\}$$ is a set of full measure, as is
$\cap_{n \geq 0} \hat{F}_A^{-n}(\G)$, since $\hat{F_A}$ is measure
preserving.  Pick $\underline{y} \in \cap_{n \geq 0} \hat{F}_A^{-n}(\G)$ and
let $\hat{g}^{-1}(\underline{y}) = \{\eta_1, \ldots, \eta_m\}$.  Then
$\Phi^n(\eta_i)$ and $\Phi^n(\eta_j)$ are nowhere coincident for all $n
\geq 0$ and $i \neq j$.  Choose  $\{t_k\}$ so that $\eta_1 - t_k \omega_R
\arrow \gamma := \gamma_1$, and $\eta_i - t_k \omega_R \arrow \gamma_i \in
\T_\vp^S$ for each $i \in \{2, \ldots, m\}$.  Then $\gamma_1,
\ldots,
\gamma_m
$  are pairwise nowhere coincident  and $\{\gamma_1,
\ldots,
\gamma_m, \gamma\}
\subset \hat{g}^{-1}(\hat{g}(\gamma))$.  Using minimality again, choose
$s_k \arrow \infty$ so that $\gamma_i - s_k \omega_R \arrow \eta_{j(i)}$ and 
$\gamma' - s_k \omega_R \arrow \eta'$ with $i \mapsto j(i)$ a bijection on
$\{1, \ldots , m\}$.  Then $\eta' = \eta_{j(i)}$ for some $i_0$, so that
$\gamma'$ is proximal to $\gamma_{i_0}$.  Either $\gamma_{i_0}$
coincides with $\gamma'$ along the segment $I$, in which case $\gamma_{i_0}$
is coincident with $\gamma_1$ along $I$, so that  $\gamma_{i_0} = \gamma_1 =
\gamma$ and we have $\gamma' \sim_p \gamma$, or the segment $I$ occurs in a
bubble formed by $\gamma_{i_0}$ and $\gamma'$.  By (2), this bubble must
terminate with coincidence.  That is, there is $t \in \R$ so that $E^s + t
\omega_R$ meets
$I$ in its interior, and there is $n \in \N$ so that $\Phi^n(\gamma_{i_0} - t
\omega_R)$ and $\Phi^n(\gamma' - t
\omega_R)$ are coincident along a common strand meeting $E^s$.  But so are 
 $\Phi^n(\gamma' - t
\omega_R)$  and  $\Phi^n(\gamma - t
\omega_R)$.  Thus $\Phi^n(\gamma_{i_0})$ and $\Phi^n(\gamma_{1}) =
\Phi^n(\gamma)$ are coincident along a common strand $J$ meeting $E^s +
\lambda^n t \omega_R$.  Now $\Phi^n(\eta_{i_0} - t_k
\omega_R)$ and $\Phi^n(\eta_1 - t_k
\omega_R)$ converge to $\Phi^n(\gamma_{i_0})$ and $\Phi^n(\gamma)$
respectively,  so for large $k$, $\Phi^n(\eta_{i_0} - t_k
\omega_R)$ and $\Phi^n(\eta_1 - t_k
\omega_R)$ must coincide along $J$ as well.  Thus $\Phi^n(\eta_{i_0})$ and
$\Phi^n(\eta_1 )$  coincide along $J + \lambda^n t_k \omega_R$, so that $i =
1$ and $\gamma \sim_p \gamma'$.
\qed

 \  

The next result for the case in which $\vp$ is unimodular is Theorem 
6.1 of \cite{mjarek}.  We need a generalization to the
nonunimodular case.  Because the proof is somewhat technical, we include it
in an appendix.

\begin{theorem}(Asubharmonicity): \label{asubhar} Suppose that $\vp$ is
strong Pisot with abelianization $A$.  If $B$ and $R$ are nonsingular integer
matrices with $AR = RB$, $p$ is a continuous surjection so that the
diagram 

{\parindent=0pt
\begin{picture}(25, 80)(-130, -40)

\put(80, -8){$\hat{F}_R$}
\put(23, -32){$\inv F_A$}
\put(-10, -8){$\hat{g}$}

\put(-25, 20){$\T^S_\vp$}
\put(75, 20){$\inv F_B$}
\put(-10,10){\vector(1,-1){30}}
\put(5,25){\vector(1,0){60}}
\put(25, 30){$p$}

\put(85,10){\vector(-1,-1){30}}

\end{picture}

commutes,} and $p, \hat{F}_R$, and $\hat{g}$ semi-conjugate the $\R$- and
$\Z$-actions on the various spaces, then $\hat{F}_R$ is a homeomorphism.

\end{theorem}

\begin{theorem}\label{closedGCC} Suppose that $\vp$ is strong Pisot.  Then $\vp$ 
satisfies GCC if and only if proximality is a closed  relation  on $\T_\vp$,
in which case $\T_\vp/\sim_p\ \simeq \inv F_a$.
\end{theorem}

Proof: If $\vp$ satisfies GCC, then $\sim_p$ is a closed equivalence relation
whose equivalence classes  are precisely the preimages of
points under the geometric realization $\hat{g}$ (Corollary
\ref{proxequiv}), so that 
$\T_\vp/\sim_p
\ \
\simeq
\  \inv F_a$.

Suppose that $\sim_p$ is a closed relation.  According to Lemma
\ref{closed}, $\sim_p$ is an equivalence relation. 
Using the strand  space
model of $\T_\vp$, $\T^S_\vp$, we have a commuting diagram of continuous 
surjections

{\parindent=0pt
\begin{picture}(25, 80)(-130, -40)

\put(80, -8){$c$}
\put(23, -32){$\inv F_A$}
\put(-10, -8){$\hat{g}$}

\put(-50, 20){$\T_\vp \simeq \T^S_\vp$}
\put(70, 20){$\T^S_\vp/\sim_p $}
\put(-10,10){\vector(1,-1){30}}
\put(-0,25){\vector(1,0){65}}
\put(35, 30){$p$}

\put(85,10){\vector(-1,-1){30}}

\end{picture}

in which} $p$ is the  quotient map $\gamma \arrow [\gamma] := \{
\gamma' :
\gamma' \sim_p \gamma\}$,
$\hat {g}$ is geometric realization, and $c([\gamma]) := \hat{g}(\gamma)$ is 
well-defined by Proposition \ref{proxgeom}.  
Note that the maps in the
diagram also  commute with the
$\Z$- and $\R$-actions on the various spaces (which,  on $\T_\vp^S/\sim_p$,
are $[\gamma] \mapsto [\Phi(\gamma)] := \Phi[\gamma]$ and $ ([\gamma],
t)
\mapsto [\gamma]_t := [\gamma - t\omega_R]$,  respectively).

{\em Claim:  }  $c$ is a covering map.

  According to Theorem \ref{GCCeq},
$\hat{g}$, and hence $c$, is bounded-to-one.   Choose
$\underline{x} \in \inv F_A$ with $m = \# c^{-1} (\underline{x}) < \infty$
maximal.  Let $c^{-1} (\underline{x}) = \{ [\gamma_1], \ldots, [\gamma_m]\}$
and choose $\underline{x}' \in \inv F_A$. Since the flow  on
$
\inv F_A$ is minimal,  there is
$\{t_n\} \subset \R$ with 
$(\underline{x})_{t_n}\arrow  \underline {x}'$. For some 
 subsequence (which without loss of  generality, and to simplify
notation, we assume to be $\{t_n\}$  itself),  $\{\gamma_i -
t_n\omega_R\}$ converges for each $i \in \{1, \ldots, m\}$, say  to
$\gamma'_i \in \T_\vp^S$.  Since  $\gamma_i$ and $ \gamma_j$
are not proximal for $i \neq j$,  there is $\delta > 0$ so  that
$d(\gamma_i - t\omega_R,
\gamma_j - t\omega_R) \geq \delta$ for  all $t \in \R$ and for all $i \neq j
\in \{1,
\ldots,  m\}$.  It  follows that $d(\gamma'_i - t\omega_R, \gamma'_j -
t\omega_R) \geq  \delta$ for all $t \in \R$ and  $i \neq j$, so that
$[\gamma'_i] \neq  [\gamma'_j]$ for $i \neq j$.  Thus $\#
c^{-1}(\underline{x}')  \geq m$. By maximality of $m$, $\#
c^{-1}(\underline{x}') = m$ and $c$ is $m$-to-one  everywhere.

Suppose that for some $\gamma \in \T^S_\vp$, $c$ is not one-to-one on any
neighborhood of $[\gamma]$.  There are then $\{[\gamma_n]\}$,
$\{[\gamma'_n]\}$ converging to $[\gamma]$ with $[\gamma_n] \neq [\gamma'_n]$
and $c([\gamma_n]) = c([\gamma'_n])$ for all $n \in \N$.  Without loss of
generality, $\gamma_n \arrow \eta$ and $\gamma'_n \arrow \eta'$.  Since
$\sim_p$ is closed, $\eta, \eta' \in [\gamma]$; i.e.,  $\eta, \eta'$ are
proximal.  According to (4) of Lemma
\ref{closed}, there is $\delta > 0$ so that $d(\gamma_n - t \omega_R,
\gamma'_n - t \omega_R) \geq \delta$ for all $t \in \R$ and all $n \in
\N$.  On the other hand, if $t \in \R$ is chosen so that  $d(\eta - t
\omega_R,
\eta' - t \omega_R) < \delta$, then for large enough $n$,  $d(\gamma_n - t
\omega_R,
\gamma'_n - t \omega_R) < \delta$.  This contradiction proves that $c$ is
locally one-to-one and hence  an
$m$-to-one covering map.

By Lemma \ref{dbyd}, there are $d \times d$ nonsingular integer  matrices $B,
R$ with $RB = AR$ and a homeomorphism $h$ so that the  diagram

{\parindent=0pt
\begin{picture}(25, 80)(-130, -40)

\put(80, -8){$\hat{F}_R$}
\put(23, -32){$\inv F_A$}
\put(-10, -8){$c$}

\put(-50, 20){$ \T^S_\vp/\sim_p$}
\put(70, 20){$\inv F_B $}
\put(-10,10){\vector(1,-1){30}}
\put(-0,25){\vector(1,0){65}}
\put(35, 30){$h$}

\put(85,10){\vector(-1,-1){30}}

\end{picture}

commutes.}   As $B = R^{-1}AR$,  $B$
is also Pisot: Let $\omega_{R, B}$ be the right eigenvector of $B$  with $R
(\omega_{R, B}) =
\omega_{R, A}$, the right Perron-Frobenius eigenvector of $A $.  Note that
$h$ and $\hat{F}_R$ semi-conjugate the flows on $\T_ \vp^S/\sim_p$ and $\inv
F_A$ with $(\underline{y}, t) \mapsto (y_1 - t  \omega_{R, B}, y_2 -
(t/\lambda)
\omega_{R, B}, \ldots) = (\underline{y}) _t$ on $\inv F_B$, and that $h$
semi-conjugates the action $[\gamma] \mapsto [\Phi(\gamma)]$ on $
\T^S_\vp/\sim_p$ with $\hat{F}_B$ on $\inv F_B$.  Now we have the commuting
diagram

{\parindent=0pt
\begin{picture}(25, 80)(-130, -40)

\put(80, -8){$\hat{F}_R$}
\put(23, -32){$\inv F_A$}
\put(-10, -8){$\hat{g}$}

\put(-30, 20){$ \T^S_\vp$}
\put(75, 20){$\inv F_B $}
\put(-10,10){\vector(1,-1){30}}
\put(-0,25){\vector(1,0){65}}
\put(25, 30){$h\circ p$}

\put(85,10){\vector(-1,-1){30}}

\end{picture}

of} surjections that semi-conjugate all of the $\Z$- and $\R$- actions.  By
the `asubharmonicity theorem' (Theorem 6.1 of \cite{mjarek} in the
unimodular case, and Theorem \ref{asubhar} of this paper for  the general
case),
$\hat{F}_R$ is one-to-one.  Thus $c$ is a  homeomorphism and $\T_\vp/\sim_p
\ 
\simeq \inv F_A$.  

 If  $\hat{g}(\gamma) =
\hat{g}(\gamma')$, then $c(p(\gamma)) =
c(p(\gamma'))$, so that  $p(\gamma) =
p(\gamma')$.  That is, $\gamma$ and $\gamma'$ are proximal, hence 
 $\gamma $ and $ \gamma'$  share an edge.  Then  $cr_\vp = 1$, and, by Theorem
\ref{GCCeq},  $\vp$ satisfies GCC.
\qed

\begin{example}Proximality not closed
\end{example}
As the Morse-Thue example ($\vp(1) = 12$, $\vp(2) = 
21$) shows, proximality need not be a closed relation,
even for a weak Pisot substitution.  Also, the fixed words represent tilings
that are proximal in one direction but not the other. \qed 

  \section{Appendix}

We complete the proofs of Lemma \ref{onlypairs} and Theorem \ref{asubhar}.

To prove 

{\parindent=0pt {\bf Lemma  \ref{onlypairs} }}{\em The 
balanced pairs $\larr 1 \\ 1 \rarr$, 
$\larr 2 \\ 2
\rarr$,  $\larr 21 \\ 12
\rarr$, 
$\larr 211
\\ 122
\rarr$, and their duals are the only essential balanced pairs for $\psi$.}

{\parindent=0pt we} prove a slightly stronger statement:  

{\parindent=0pt {\bf
Lemma  \ref{onlypairs}$'$}}{\em 
\label{restate} Given any balanced pair $\larr u \\ v \rarr$ for $\psi$,
$\larr
\psi(u) \\ \psi(v) \rarr$ is the concatentation of the irreducible balanced
pairs
$\larr 1 \\ 1 \rarr$,  $\larr 2 \\ 2
\rarr$,  $\larr 21 \\ 12
\rarr$, 
$\larr 211
\\ 122
\rarr$, $\larr 21211 \\ 11212 \rarr$, and their duals.} 

 Since $\larr
\psi(21)
\\
\psi(12) 
\rarr$, $\larr
\psi(211)
\\
\psi(112) 
\rarr$ and 
$\larr
\psi(21211)
\\
\psi(11212) 
\rarr$ are the concatenation of $\larr 1 \\ 1 \rarr$,  $\larr 2 \\ 2
\rarr$,  $\larr 21 \\ 12
\rarr$, 
$\larr 211
\\ 122
\rarr$, and their duals, it follows that for any balanced pair $\larr u \\ v
\rarr$, $\larr
\psi^2(u)
\\
\psi^2(v)
\rarr$ is also such a concatentation, and Lemma \ref{onlypairs} is proved.

Proof of Lemma  \ref{onlypairs}$'$:  The following proof consists of
demonstrating that for any irreducible balanced pair $\larr u \\ v \rarr = \larr
u_1 \ldots u_n \\ v_1 \ldots v_n \rarr $, the balanced pair $\larr \psi(u) \\
\psi(v)
\rarr $ can be factored into a product of irreducible balanced pairs by first
reducing  $\larr \psi(u_1) \\
\psi(v_1)
\rarr $ (with remainder), then reducing  $\larr \psi(u_1u_2) \\
\psi(v_1v_2)
\rarr $ (with remainder), etc; each step in this sequence is represented by an
edge (and the adjoining vertices) in a finite graph $G$.  The complete
factorization of 
$\larr
\psi(u)
\\
\psi(v)
\rarr $ is then represented by a path of length $n$ in $G$ with the $i^{th}$ edge
of the path labeled $\larr u_i \\ v_i \rarr $ (and the adjoining vertices
indicating, in part, remainders in the reductions of $\larr \psi(u_1\ldots u_i) \\
\psi(v_1\ldots v_i)
\rarr $ and $\larr \psi(u_1\ldots u_{i+1}) \\
\psi(v_1\ldots v_{i+1})
\rarr $.  The proof will be completed
by noting that for any
$\larr u
\\ v
\rarr
$, only the vertices and edges of the graph $G$ will be visited, and the
irreducible balanced pairs arising in each step of the factorization are included
in those listed in Lemma
\ref{onlypairs}$'$ above.

 Given two words $u, v$ of the same length allowed for the substitution
$\psi$, we define the discrepancy vector of
$u$ and $v$ to be the difference of the content vectors of $u$ and $v$, $
l(u) - l(v)$.   The discrepancy vector has integer entries, and
since $\psi$ is on a two letter alphabet, is of the form
$(m, -m)$, where $m \in \Z$.  Let 
$dis(u, v) = |m|$ be the {\em  discrepancy number} of $u$ and
$v$.  In the case of $\psi$, $|m| = 0, 1$ or 2. Note that $dis(u, v) = 0$ if and
only if $\larr u \\ v \rarr$ is balanced.

Let $\larr u \\ v \rarr = \larr u_1u_2 \ldots u_n \\ v_1 v_2 \ldots v_n
\rarr$ be an irreducible balanced pair for
$\psi$. Without loss of generality, $u_1 = 2, v_1 = 1$.  Since
$|\psi(2)| = 7$ while
$|\psi(1)| = 19$, and $\larr u \\ v \rarr$ is irreducible, $|\psi(u_1 u_2
\ldots u_i)|  < |\psi(v_1 v_2 \ldots v_i)| $ for all $1 \leq i < n$, and
$dis(u_1 \ldots u_i, v_1 \ldots v_i)$ indicates how many
more 2's appear in $u_1 \ldots u_i$  than in $v_1 \ldots v_i$.

In the graph $G$ below, each edge is labelled by a pair of symbols $\lsqarr a_1 \\
a_2
\rsqarr $. Each vertex has two components: 
an integer representing a discrepancy number, and a
pair of words
$\lsqarr w_1
\\ w_2
\rsqarr$ representing, under the right circumstances, the remainder in the
reduction of a  particular pair of words.  As described above, for an irreducible
balanced pair
$\larr 2u_2
\ldots u_n
\\ 1 v_2 \ldots v_n \rarr$, the reduction of $\larr \psi(2u_2 \ldots u_n) \\
\psi(1 v_2 \ldots v_n) \rarr$ is represented by a path of length
$n$ through $G$ that begins and ends at $(0, 0)$ in which the $i^{th}$ vertex
in the path indicates both $dis(u_1 \ldots u_i, v_1 \ldots v_i)$ and the
remainder upon writing 
$\larr
\psi(u_1u_2
\ldots u_i)
\\
\psi(v_1 v_2 \ldots v_i) \rarr$ in irreducible balanced pairs.  The edge to the
$(i+1)^{st}$ vertex is labelled with $\lsqarr a_1 \\ a_2 \rsqarr = \lsqarr u_{i +
1}
\\ v_{i + 1}\rsqarr$, and the $(i+1)^{st}$ vertex consists of $dis(u_1 \ldots
u_{i+1}, v_1 \ldots v_{i+1})$ and the remainder in the reduction of $\larr
\psi(u_1u_2
\ldots u_{i + 1})
\\
\psi(v_1 v_2 \ldots v_{i + 1}) \rarr$.  

To simplify the writing of remainders,
we use the notation
$x_i^+$ ($y_i^+$, respectively) to denote the subword  $x_i \ldots
x_{19}$ of $\psi(1)$ ($y_i
\ldots y_7$ of $\psi(2)$, respectively).
In the following diagram,  $a, b, c$ and $ d$ denote $\lsqarr 1 \\
1 \rsqarr$, $\lsqarr 1 \\
2 \rsqarr$, $\lsqarr 2 \\
1 \rsqarr$, and $\lsqarr 2 \\
2 \rsqarr$,  respectively.

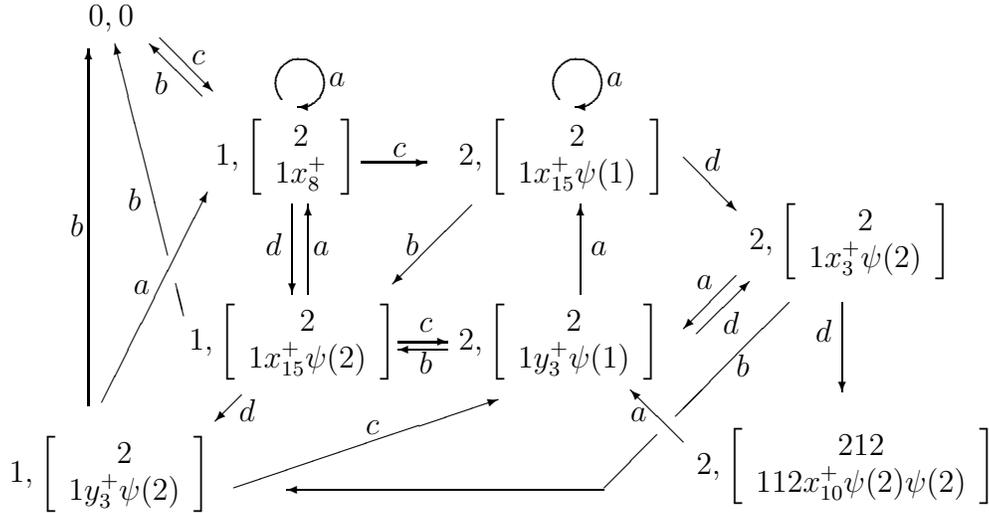
\begin{picture}(25, 230)(-80, -120)

\put(2, -7){$d$}
\put(20, -7){$a$}
\put(55, -7){$b$}
\put(50, 30){$c$}
\put(60, -36){$c$}

\put(-26, 65){$c$}
\put(125, -7){$a$}

\put(-65, 80){$0, 0$}
\put(-95, -93){$1, \lsqarr 2 \\ 1y_3^+ \psi(2) \rsqarr$}
\put(-65, -64){\vector(0, 1){135}}
\put(-72, 0){$b$}

\put(-17, 27){$1, \lsqarr 2 \\ 1x_8^+  \rsqarr$}
\put(15,58){\arc(2, -9){300}}
\put(16, 49){\vector(-1, 0){2}}
\put(26, 57){$a$}
\put(75, 27){$2, \lsqarr 2 \\ 1x_{15}^+ \psi(1) \rsqarr$}
\put(120,58){\arc(2, -9){300}}
\put(121, 49){\vector(-1, 0){2}}
\put(131, 57){$a$}

\put(185,-5){$2, \lsqarr 2 \\ 1x_{3}^+ \psi(2) \rsqarr$}
\put(165,-90){$2, \lsqarr 212 \\ 112x_{10}^+ \psi(2) \psi(2) \rsqarr$}
\put(-27, -43){$1, \lsqarr 2 \\ 1x_{15}^+ \psi(2) \rsqarr$}
\put(-35, -7){\vector(-1, 4){20}}
\put(-29, -30){\line(-1, 4){3}}
\put(-50, 10){$b$}
\put(75, -43){$2, \lsqarr 2 \\ 1y_3^+ \psi(1) \rsqarr$}

\put(12,12){\vector(0, -1){34}}
\put(18,-22){\vector(0, 1){34}}
\put(220, -25){\vector(0, -1){34}}
\put(121, -22){\vector(0, 1){34}}
\put(52,-40){\vector(1,0){19}}

\put(38,28){\vector(1,0){25}}

\put(80, 12){\vector(-1, -1){30}}
\put(71, -43){\vector(-1, 0){19}}
\put(60, -52){$b$}
\put(-10, -95){\vector(3, 1){100}}
\put(40, -75){$c$}
\put(-60, -63){\vector(1, 2){40}}
\put(-48, -22){$a$}
\put(-38, 75){\vector(1, -1){20}}
\put(-22, 53){\vector(-1, 1){20}}
\put(-40, 54){$b$}
\put(-7, -60){\vector(-1, -1){10}}
\put(-8, -70){$d$}
\put(180, -15){\vector(-1, -1){20}}
\put(210, -40){$d$}
\put(165,-37){\vector(1, 1){20}}
\put(165, -20){$a$}
\put(175,-37){$d$}
\put(160, 30){\vector(1, -1){20}}
\put(168, 25){$d$}
\put(160, -78){\vector(-1, 1){20}}
\put(140, -70){$a$}
\put(200, -25){\line(-1, -1){43}}
\put(150, -75){\line(-1, -1){20}}
\put(130, -96){\vector(-1, 0){120}}
\put(180,-53){$b$}

\end{picture}

Some transitions between vertices are not allowed.  For instance,
the edge $c = \lsqarr 2 \\ 1\rsqarr$
cannot leave a vertex with discrepancy  2, since this would result in a
discrepancy of 3, not allowed for $\psi$.  Also,  the edges
$b  =
\lsqarr 1
\\ 2 \rsqarr$ and 
$ d  = \lsqarr 2 \\
2 \rsqarr$ cannot leave the vertices with terminal factors of $\psi(2)\psi(2)$ or 
$y_3^+\psi(2)$ in their remainders, since the subword 222 is not allowed for
$\psi$.  Since every other edge appears in $G$, every irreducible balanced pair
is represented by a path in $G$.

We illustrate the use of $G$ with the balanced pair $\larr 211
\\ 112 \rarr$.  The path representing the reduction of $\larr \psi(211)
\\ \psi(112) \rarr$ will involve edges
labelled by
$c =
\lsqarr 2
\\ 1
\rsqarr$,
$a = \lsqarr 1 \\ 1 \rsqarr$, and
$ b = \lsqarr 1
\\ 2
\rsqarr$, in that order, discrepancy numbers of 1, 1, and 0, and
remainders of $ \lsqarr 2 \\ 1x_8^+  \rsqarr$, $ \lsqarr 2 \\ 1x_8^+ 
\rsqarr$, and $0$.

We invite the reader to verify some of the calculations. Note that movement
from one vertex to another can be checked independently of the remainder of
the graph.  

We leave the reader to check  that all
reductions except  that involving the single edge from $2, \lsqarr 2 \\
112v_{10}^+ \psi(2) \psi(2) \rsqarr$ involve only the irreducible balanced pairs
described in the statement of Lemma \ref{onlypairs}; this exception also
includes $\larr 21211 \\ 11212 \rarr$ in its reduction.  It is
 easy to verify that the reduction of $\larr \psi(21211) \\ \psi(11212) \rarr$
also involves only the balanced pairs of Lemma \ref{onlypairs}.  This proves the
lemma.
\qed

\

We now prove Theorem \ref{asubhar}.

\

{\parindent=0pt {\bf Theorem \ref{asubhar}}{\em (Asubharmonicity): Suppose that
$\vp$ is strong Pisot with abelianization $A$.  If $B$ and $R$ are nonsingular
integer matrices with $AR = RB$, $p$ is a continuous surjection so that the
diagram 

{\parindent=0pt
\begin{picture}(25, 80)(-130, -40)

\put(80, -8){$\hat{F}_R$}
\put(23, -32){$\inv F_A$}
\put(-10, -8){$\hat{g}$}

\put(-25, 20){$\T^S_\vp$}
\put(75, 20){$\inv F_B$}
\put(-10,10){\vector(1,-1){30}}
\put(5,25){\vector(1,0){60}}
\put(25, 30){$p$}

\put(85,10){\vector(-1,-1){30}}

\end{picture}

commutes, } and $p, \hat{F}_R$, and $\hat{g}$ semi-conjugate the $\R$- and
$\Z$-actions on the various spaces, then $\hat{F}_R$ is a homeomorphism.}
}

\

Proof:  Given a strand (finite or infinite) $\gamma$ that meets $E^s$, let
$\hat{\gamma}$ denote the state determined by $\gamma$ (i.e., the edge of
$\gamma$, closed on the initial end and open on the terminal end, that meets
$E^s$).  We write $\gamma \sim_0 \eta$ if $\gamma $ and $\eta$ are strands
that meet
$E^s$ for which  $\widehat{\Phi^n(\gamma)} = \widehat{\Phi^n(\eta)}$ for some
$n
\in N$.  Note that $\widehat{\Phi^(\hat{\gamma})} =
\widehat{\Phi^(\gamma)}  $. Given $\gamma \in \T^S_\vp$, let $\F_\gamma
= \{ v \in \cup_{n \geq 0} A^{-n}\Z^d : \gamma \sim_0 \gamma + v\}$. 

Recall that $\omega_R$ is the right eigenvector for $A$.  Define $$W^s(\gamma)
= \{
\eta
\in \T^S_\vp : d(\Phi^n(\gamma), \Phi^n(\eta)) \arrow 0\ \textrm{as} \ n
\arrow
\infty\}$$ and let $pr_A$ denote the projection of $\R^d$ onto $E_A^u$ along
$E^s_A$. Letting $$ret(\gamma) = \{ t \in \R : \gamma - tw_R \in
W^s(\gamma)\}$$ be the set of return times for $\gamma$, we see that, as long
as
$\Phi^n(\gamma)$ does not have a vertex on $E^s$ for any $n \geq 0$, $v \in
\F_\gamma $ if and only if $pr_Av = t\omega_R$ for some $t \in ret(\gamma)$. 
Let $$\G = \{ \gamma \in \T^s_\vp: \Phi^n(\gamma)  \textrm{ does not have
a vertex on } E^s  \textrm{ for any } n \geq 0\}.$$  Then $\G$ has full
measure in $ \T_\vp^s$.  Let $H = \langle \cup_{\gamma \in \G}
\F_\gamma \rangle$ be the subgroup of $\cup_{n
\geq 0} A^{-n}\Z^d $ generated by $\cup_{\gamma \in \G} \F_\gamma$.

Since $(\Phi^{n+1}(\Phi^{-1}(\gamma) + A^{-1}v)) \hspace*{.02in} \widehat{  
} = (\Phi^n(\gamma + v)) \hspace*{.02in} \widehat{   } $, and $\gamma \in
\G$ if and only if
$\Phi^{-1}(\gamma) \in \G$, we can conclude that:

\

{\parindent=0pt (A.1)  $A^{-1}(H) \subset H$.}

\

For $i \in \{1, \ldots, d\}$, let 
\begin{eqnarray*}
\Theta(i) & = \{v \in \Z^d: \textrm{ for
some (and hence any) }  \gamma \in \T^S_\vp,  \textrm{ there are edges } \\ &
I, J
\textrm{ of } \gamma, \textrm{ both of type } i, \min(J) - \min(I) = v\} 
\end{eqnarray*}
be the set of return vectors for type $i$.  Note that for $v \in \Theta(i)$
and
$\gamma \in \T_\vp^S$, $\gamma$ and $\gamma + v$ share an edge, say $I$.  It
follows that if $E^s + t \omega_R$ meets $I$ in its interior, and $\gamma -
t \omega_R \in \G$ (as will be the case for most such $t$), then $(\gamma -t
\omega_R) \sim_0 (\gamma -t
\omega_R) + v$, hence $v \in \F_{\gamma -t
\omega_R}$.  Thus:

\

{\parindent=0pt (A.2)  $\Theta(i) \subset H$ for} all  $i \in \{1, \ldots,
d\}$.

\

Now fix $\gamma \in \G$, and define $$[v + H]^+ = \{i: i \textrm{ is the
type of an edge } I \textrm{ of }  \gamma \textrm{ with} \min I \in v +
H\},$$
$$[v + H]^- = \{i: i \textrm{ is the type of an edge } I \textrm{ of }
\gamma \textrm{ with} \max I \in v + H\}.$$  If $i \in [v + H]^+
\cap [u + H]^+$, then $\min I + H = v + H$ and $\min J + H = u + H$ for some
edges $I, J$ of $\gamma$ of type $i$.  Then $u - v + H = \min J - \min I + H
= H$, since $\min J - \min I \in \Theta(i) \subset H$, hence $u - v \in
H$.  A similar statement can be made if $i \in   [v + H]^-
\cap [u + H]^-$, so that:

\

{\parindent=0pt
(A.3) } If $[v + H]^+ \cap [u + H]^+ \neq \emptyset$ or $[v + H]^- \cap [u +
H]^- \neq \emptyset$, then $u - v \in H$.

\

Recall that the flow on $\T^S_\vp$ is uniquely ergodic, and that if for each 
$i
\in \{1, \ldots, d\}$,  $f_i$ denotes the frequency of occurrence  of tiles
of type
$i$ in
$\gamma \in \T^S_\vp$, $$f_i := \lim_{N \arrow \infty} \frac{1}{2N} \mu \{t
\in
[-N, N]: (\gamma - t\omega_R)\widehat{} \textrm{ has type } i\},$$ then
$(f_1,
\ldots, f_d)^{tr}$ is a right Perron-Frobenius eigenvector for $A$.
Since the characteristic polynomial of $A$ is irreducible over $\Q$, the
entries of $(f_1, \ldots, f_d)$ are independent over $\Z$.  It follows from
(A.2) that an edge of type $[v + H]^- $ in $\gamma$ must be followed by an
edge of type $[v + H]^+$, and an edge of type $[v + H]^+$ must be preceded
by an edge of type $[v + H]^-$.  Thus occurrences in $\gamma$ of edges of
type $[v + H]^-$ are in one-to-one correspondence with occurrences of edges
of type
$[v + H]^+$.  That is,  the frequency of type $[v + H]^+$ equals the
frequency of type $[v + H]^-$:  $$\Sigma_{i \in [v + H]^+} f_i = \Sigma_{i
\in [v + H]^-} f_i .$$  As the $f_i$ are linearly independent over $\Z$, it
must be that:

\

{\parindent=0pt (A.4) } $[v + H]^+ = [v + H]^-$.

\

Now fix  $i \in \{1, \ldots, d\}$, and let $e_i = (0, \ldots, 1, \ldots, 0)$
be the unit vector in the $i$ direction.  Let $I$ be an edge of $\gamma$ of
type $i$, $u = \min I$ and $v = \max I$.  Then $i \in [u + H]^+ = [u + H]^-$
and $i \in [v + H]^-$.  Thus $v - u \in H$ by (A.3) so that $e_i \in H$ and
thus $\Z^d \subset H$.  From (A.1) we now have $\cup_{n \geq 0} A^{-n}\Z^d
\subset H$.  Thus:

\

{\parindent=0pt
(A.5) } $H = \langle \cup_{n \geq 0} A^{-n}\Z^d \rangle$.

\

Suppose that $t \in ret(\gamma)$ (i.e., $t$ is a return time for $\gamma $.) 
Then $t$ is also a return time for $p(\gamma) \in \inv F_B$, since $p$
semi-conjugates the $\R$- and $\Z$-actions on $\T^S_\vp$ with those on $\inv
F_B$.  That is, if  $t \in ret(\gamma)$, then $t\omega_{R, B} \in
pr_B(\cup_{n
\geq 0} B^{-n}\Z^d) $, where $pr_B$ is the projection of $\R^d$ onto $E^u_B$
along $E^s_B$ and $\omega_{R, B} := R^{-1}(\omega_{R, A})$ is a right
eigenvector of $B$ corresponding to the Perron-Frobenius eigenvalue $\lambda
$ of $A$.

Choose $\gamma \in \G$ and $v \in \F_\gamma$.  Then for $t \in ret(\gamma)$
and $u \in \cup_{n \geq 0} B^{-n}\Z^d$, 
$$pr_Av = t\omega_{R, A} = Rt\omega_{R, B} = R pr_B(u).$$  But $R pr_B(u) =
pr_A(Ru)$ since $R^{-1}AR = B$ and, by irreducibility of the characteristic
polynomial of $A$, $v = Ru$.  That is, $\F_\gamma \subset R( \cup_{n \geq 0}
B^{-n}\Z^d)$.  Thus $$H = \langle \cup_{\gamma \in \G}  \F_\gamma \rangle
\subset R( \cup_{n \geq 0} B^{-n}\Z^d). $$  In addition, since $RB = AR$, it
follows that $ R(  B^{-n}\Z^d) \subset A^{-n} \Z^d$, so that by (A.5):

\

{\parindent=0pt (A.6)} $ \cup_{n \geq 0} A^{-n}\Z^d =  R( \cup_{n \geq 0}
B^{-n}\Z^d)$.

\

Since $\hat{F}_R $ is a covering map and $\inv F_B$ is compact,
$\hat{F}_R^{-1}(\underline{0} := (0, 0, \ldots))$ is finite.  Let
$\underline{x}  \in \hat{F}_R^{-1}(\underline{0} )$.  As  $\hat{F}_R $ 
semi-conjugates  $\hat{F}_B $ with  $\hat{F}_A $, $\underline{x} $ must be
periodic under $\hat{F}_B $; without loss of generality, assume
$\underline{x} = (x, x, \ldots )$ is fixed, say $x = y + \Z^d \in \cT^d$. 
Since $\hat{F}_R(\underline{x}) = \underline{0}$, $Ry = v \in \Z^d$.  Also,
$\underline{x} \in \inv F_B$, hence $B(y) = y + w$, where  $w \in \Z^d$.
We have 

{\parindent=0pt
\begin{picture}(25, 95)(-140, -50)

\put(0, -10){$R$}
\put(75, -10){$R$}
\put(40, 25){$B$}
\put(40, -25){$A$}

\put(-15, 17){$y + w$}
\put(75, 17){$y$}
\put(-5, -33){$Av$}
\put(75, -33){$v$}


\put(15,12){\vector(0, -1){37}}
\put(65,-30){\vector(-1,0){45}}
\put(70,12){\vector(0,-1){37}}
\put(65,20){\vector(-1,0){45}}

\put(13,12){\line(1,0){4}}
\put(65,-32){\line(0,1){4}}
\put(68,12){\line(1,0){4}}
\put(65,18){\line(0,1){4}}

\end{picture}

That is}, $Av = R(y + w)$.  Since $R(y + w) \in \Z^d$,  $v = A^{-1}( R(y
+ w))
\in
\cup_{n
\geq 0} A^{-n}
\Z^d$.  Thus, by (A.6), there are  $ n \in \N$ and $u \in \Z^d$ so that 
$RB^{-n} u = v$.
Now $y = R^{-1}v = B^{-n}u$, so that
 $$By = B^{-n+1}u = y + w,$$  $$y = B^{-n+1}u - w,$$ $$By =
B^{-n+2}u - Bw = y + w,$$  and $$y = B^{-n+2}u - Bw - w.$$
Continuing, we obtain 
$$y = u - B^{n-2}u - B^{n-3}u - \ldots - Bw - w \in \Z^d. $$ That is, 
$\underline{x} = \underline{0}$.  Hence $\# \hat{F}^{-1}_R(0) = 1$.  
 Since $\hat{F}_R$ is a covering map, and $\inv F_B$ is
connected, $\hat{F}_R$ is a homeomorphism.

 \qed

\

{\small {\parindent=0pt Department of Mathematics, Montana State University,
Bozeman, MT 59717

barge@math.montana.edu

\

Department of Mathematics, College of Charleston, Charleston, SC 29424

diamondb@cofc.edu}

}

\end{document}